\documentclass{amsart}
\usepackage{amssymb}

\usepackage{xypic}
\xyoption{all}

\newtheorem{thm}{Theorem}[section]
\newtheorem*{thm*}{Theorem}
\newtheorem{prop}[thm]{Proposition}
\newtheorem{lem}[thm]{Lemma}
\newtheorem{lem-def}[thm]{Lemma-Definition}

\newtheorem{cor}[thm]{Corollary}
\newtheorem*{cor*}{Corollary}
\newtheorem*{conj*}{Conjecture}

\theoremstyle{definition}
\newtheorem{defn}[thm]{Definition}
\newtheorem{example}[thm]{Example}

\newtheorem{rem}[thm]{Remark}
\newtheorem{remark}[thm]{Remark}

\newcommand{\op}{\operatorname}
\renewcommand{\gg}{{\op{\mathfrak g}}}

\newcommand{\sdp}{\mathbin{{>}\!{\triangleleft}}}

\newcommand{\bbA}{\mathbb{A}}
\newcommand{\bbG}{\mathbb{G}}

\newcommand{\NC}{{\op NC}}
 \newcommand{\Ker}{{\op{Ker}}}

\newcommand{\ZZ}{{\mathbb Z}}

\newcommand{\CC}{{\mathbb C}}

\renewcommand{\AA}{{\mathbb A}}
\newcommand{\GG}{{\mathbb{G}}}

\newcommand{\be}{\begin{enumerate}}
\newcommand{\ee}{\end{enumerate}}

\newcommand{\id}{\op{id}}

\newcommand{\SL}{\op{SL}}
\newcommand{\GL}{\op{GL}}

\newcommand{\RT}{\op{RT}}

\newcommand{\OO}{\op{O}}

\newcommand{\PGLn}{\op{PGL}_n}

\newcommand{\tr}{\op{tr}}

\newcommand{\Aut}{\op{Aut}}

\newcommand{\diag}{\op{diag}}

\newcommand{\Ad}{\op{Ad}}
\newcommand{\Rad}{\op{Rad}}

\newcommand{\Lie}{\op{Lie}}

\newcommand{\lra}{\to}

\newcommand{\Slice}{\op{Slice}}
\newcommand{\Sym}{\op{S}}

\newcommand{\Stab}{\op{Stab}}

\newcommand{\Mn}{{\op{M}}_n}
\newcommand{\M}{{\op{M}}}
\newcommand{\Mat}{{\op{M}}}

\DeclareMathOperator{\catq}{/\!/}

\title{Is the Luna stratification intrinsic?}
\author[J. Kuttler and Z. Reichstein]{J. Kuttler and Z. Reichstein}
\address{Department of Mathematics, University of British Columbia,
Vancouver, BC V6T 1Z2, Canada\\ 
{\em Current address:} Department of Mathematical and Statistical Sciences, University of Alberta, Edmonton, AB T6G 2G1, Canada}
\email{jochen.kuttler@ualberta.ca}
\address{Department of Mathematics, University of British Columbia,
Vancouver, BC V6T 1Z2, Canada}
\thanks{J. Kuttler was a PIMS Postdoctoral Fellow at UBC}
\thanks{Z. Reichstein was partially supported by an NSERC research grant}
\email{reichst@math.ubc.ca}


\begin{document}

\begin{abstract} Let
$G \lra \GL(V)$ will be a finite-dimensional linear representation of
a reductive linear algebraic group $G$ on a finite-dimensional 
vector space $V$, defined over an algebraically closed field 
of characteristic zero. The categorical quotient $V \catq G$
carries a natural stratification, due to D. Luna.
This paper addresses the following questions:

\smallskip
(i) Is the Luna stratification of $X$ intrinsic? That is,
does every automorphism of $V \catq G$ map each stratum 
to another stratum? 

\smallskip
(ii) Are the individual Luna strata in $X$ intrinsic?
That is, does every automorphism of $V \catq G$ maps each stratum
to itself?

\smallskip
In general, the Luna stratification is not intrinsic.
Nevertheless, we give positive answers to questions (i) and (ii)
for large classes of interesting representations. 
\end{abstract}

\maketitle

% \tableofcontents

\section{Introduction}

Throughout this paper $k$ will be an algebraically closed field 
of characteristic zero, $G \lra \GL(V)$ will be a representation of
a reductive linear algebraic group $G$ on a finite-dimensional 
vector space $V$, defined over $k$, and
$\pi \colon V \lra X = V \catq G$ will denote the categorical quotient map
for the $G$-action on $V$. For the definition and a discussion of the
properties of the categorical quotient in this setting, see, 
e.g.,~\cite{mf}, \cite{pv} or \cite{kraft}.

There is a natural stratification on $X$, due to
D. Luna; we shall refer to it as {\em the Luna stratification}.
Recall that for every $p \in X$ the fiber $\pi^{-1}(p)$ has a unique
closed orbit. Choose a point $v_p$ in this orbit. Then the stabilizer
subgroup $\Stab(v_p)$, is reductive, and its conjugacy class in $G$
is independent of the choice of $v_p$. This subgroup determines
the stratum of $p$. More precisely, the Luna stratum associated
to the conjugacy class $(H)$ of a reductive subgroup $H \subseteq G$
is defined as
\[ X^{(H)} = \{ p \in X \, | \Stab(v_p) \in (H) \} \, . \]
There are only finitely many Luna strata, and each
stratum is a locally closed non-singular subvariety of $X$.
Moreover, if we set
\[V^{\langle H \rangle} = \{ v \in V \mid \; \text{$G \cdot v$ is closed
and $\Stab(v) = H$} \} \]
then $\pi$ restricts to a principal $N_G(H)/H$-bundle 
$V^{\langle H \rangle} \to X^{(H)}$.  For proofs of these assertions
see \cite[Section 6.9]{pv} or~\cite[Section I.5]{schwarz}.

\smallskip
The Luna stratification provides a systematic approach to the problem
of describing the $G$-orbits in $V$; it also plays an important role in
the study of the geometry and (if $k = \CC$) the topology of the categorical
quotient $X = V \catq G$. In this paper we shall address 
the following questions.

\smallskip
(i) Is the Luna stratification of $X$ intrinsic? In other words, is it
true that for every automorphism $\sigma \colon X \lra X$ and
every reductive subgroup $H \subset G$ there is a reductive
subgroup $H' \subset G$ such that $\sigma(X^{(H)}) = X^{(H')}$?

\smallskip
(ii) Are the Luna strata in $X$ intrinsic? Here we say that $X^{(H)}$
is intrinsic if $\sigma(X^{(H)}) = X^{(H)}$
for every automorphism $\sigma \colon X \lra X$.

\smallskip
In general, the Luna stratification is not intrinsic.
Indeed, there are many examples, where $V \catq G$ is an affine
space (cf. e.g., \cite[Section 8]{pv}) and the automorphism group of an affine
space is highly transitive, 
so that points in the same stratum can be taken
by an automorphism to points in different strata.
Moreover, even in those cases where the Luna stratification
is intrinsic, the individual strata may not be;
see Example~\ref{ex.regular}.  The purpose of this paper is to
show that one can nevertheless give positive answers to 
questions (i) and (ii) for large classes of interesting 
representations. 

We begin with the case where $G$ is a finite group. Recall that
a non-trivial $g \in \GL(V)$ is called a pseudo-reflection if
$g$ has finite order and fixes (pointwise) a hyperplane in $V$.
If $G \lra \GL(V)$ is a representation and $G$ is generated 
by elements that act as pseudo-reflections on $V$ then 
by a theorem of Chevalley and Shephard-Todd,
$V \catq G$ is an affine space. As we remarked above, in this case
the Luna stratification cannot be intrinsic. To avoid this situation, 
the theorem below assumes that $G$ contains no pseudo-reflections.
In particular, this condition is automatically satisfied for
representations $G \to \SL(V) \subset \GL(V)$.

\begin{thm} \label{thm1}
Let $V$ be a finite-dimensional $k$-vector space. Suppose
a finite subgroup $G \subset \GL(V)$ contains no pseudo-reflections.
Then for every automorphism
$\sigma$ of $X = V \catq G$ there is an automorphism
$\tau$ of $G$ such that $\sigma(X^{(H)}) = X^{(\tau(H))}$
for every subgroup $H$ of $G$.
In particular, the Luna stratification in $X$ is intrinsic.
\end{thm}

A similar result in the analytic setting is due to Prill~\cite{prill}; 
our proof in Section~\ref{sect.thm1}
is an algebraic variant of Prill's argument.
Note that under the assumptions of Theorem~\ref{thm1}
the individual Luna strata are intrinsic in many cases
(but not always; see Example~\ref{ex.regular}). 
In the course of the proof we show that 
every automorphism of $X$ can be lifted to an automorphism of $V$. 

Next we turn to representations $V$ of (possibly infinite) reductive
groups $G$. There is no direct analogue of the Chevallet-Shephard-Todd 
theorem in this setting, so we replace the condition that no
element of $G$ acts as a pseudo-reflection on $V$ by the stronger 
condition that $V$ is a multiple $W^r$ of another $G$-representation 
$G \lra \GL(W)$ for sufficiently high $r$. A general (and somewhat vague)
principle in invariant theory says that replacing a $G$-variety $Z$ 
by a power $Z^r$ often ``improves" the properties of 
the underlying action, assuming $r$ is sufficiently large.  
(For two unrelated recent results along these lines 
see \cite[Corollary, p. 1606]{lorenz} 
and \cite{popov4}.) In our setting this principle 
leads to the following positive answers to question (i).

\begin{thm} \label{thm2}
Let $G \lra \GL(W)$ be a finite-dimensional linear 
representation of a reductive algebraic group $G$. 
Then the Luna stratification in $W^r \catq G$ is intrinsic if

\smallskip
{\rm (a)} $r \ge 2 \dim(W)$, or
\smallskip

{\rm (b)} $G$ preserves a nondegenerate quadratic form on $W$ and
$r \ge \dim(W) + 1$, or
\smallskip

{\rm (c)} $W = \gg$ is the adjoint representation of $G$ and $r \geq 3$.
\end{thm}

Note that in this setting we do not know under what circumstances 
an automorphism $\sigma \colon V \catq G \to V \catq G$ can
be lifted to $V$, so our proof is indirect; see Remark~\ref{rem.cs}.  
Along the way we show that if $V = W^r$ is as in Theorem~\ref{thm2} and
$S$ is a Luna stratum in $V \catq G$ then $\overline{S}$
is singular at every point $x \in \overline{S} \setminus S$.

Our final result concerns the natural $\GL_n$ action of
on the space $V = \Mn^r$ of $r$-tuples of $n \times n$-matrices
by simultaneous conjugation.
The variety $X = \Mn^r \catq \GL_n$ has been extensively studied 
in the context of both invariant and PI theories;
an overview of this research area can be found in~\cite{formanek},
\cite{procesi} or~\cite{df}. In~\cite{reichstein1, reichstein2} 
the second author constructed 
a large family of automorphisms of $X = V^r \catq \PGLn$ (for $r \ge n + 1$).
Every automorphism in that family preserved the Luna strata, so
it is natural to conjecture that the same should be true for every 
automorphism of $X$.  The following result proves this 
conjecture for any $r \ge 3$.

\begin{thm}\label{thm3}
Suppose $r \ge 3$. Then every Luna stratum in 
$X = \Mn^r \catq \GL_n$ is intrinsic.
\end{thm}

\noindent
Note that Theorem~\ref{thm3} fails if $n = 1$ or
$(n, r) = (2, 2)$; see Remark~\ref{rem.(2,2)}. 

The fact that the principal stratum $X^{(\{ e \} )}$
is intrinsic is an immediate consequence of
a theorem of LeBruyn and Procesi~\cite[Theorem II.3.4]{lbp},
which says that $X^{(\{ e \})}$ is precisely the smooth
locus of $X$. This result served both as a motivation
for and as a starting point of our proofs of Theorems~\ref{thm2} 
and~\ref{thm3}.

\section{Proof of Theorem~\ref{thm1}}
\label{sect.thm1}

In this section we will establish Theorem~\ref{thm1} 
by proving the following stronger assertion.

\begin{prop} \label{prop3.2}
Let $G$ be a finite subgroup of $\GL(V)$, 
$\pi \colon V \to X$ be the categorical
quotient map for the natural $G$-action on $V$,
and $\sigma \colon X \to X$ be an automorphism. Assume that
$G$ contains no pseudo-reflections.  Then

\smallskip
{\rm(a)} $\sigma$ lifts to an automorphism $\tilde \sigma\colon V \to V$.

\smallskip
{\rm(b)} Moreover, there exists a $\tau \in \Aut(G)$ such that
$ \tilde \sigma(gv) = \tau(g) \tilde \sigma(v)$
for all $g \in G$ and all $v \in V$.
\end{prop}

\subsection{Preliminaries}

Our proof of Proposition~\ref{prop3.2} will rely on the three
simple lemmas below.

Recall that an algebraic variety $Y$ is called \emph{simply 
connected} if it admits no nontrivial finite \'etale covers.

\begin{lem}\label{lem.sc0}
Suppose $\pi \colon E \to Z$ is a finite \'etale cover. Then
any map $f \colon Y \lra Z$, with $Y$ simply connected, factors
through $\pi$.
\end{lem}

\begin{proof} Consider the finite \'etale cover
$\pi_Y \colon Y \times_Z E \to Y$.  Since $Y$ is simply connected,
it has a section $s \colon Y \to Y \times_Z E$.
Now the commutative diagram
\[ \xymatrix{ Y \times_Z E \ar@{->}[d]  
\ar@{->}[r]_(.6){\pi_Y} &
\ar@/_1pc/[l]_{\quad s} Y
 \ar@{->}[d]^f \cr
 E \ar@{->}[r]^{\pi} &  Z  ,} \]
shows that $f$ factors through $\pi$.
\end{proof}

We will need the following well known remark in the proof of Proposition~\ref{prop3.2} below:
\begin{rem}\label{rem.sc-1}
Suppose $Y$ is an irreducible smooth simply connected variety and $U$ is an open subset of $Y$. If $Y \setminus U$ has codimension $\geq 2$ in $Y$, then $U$ is simply connected.

This is an immediate application of the fact that every 
finite \'etale cover of $U$ extends to one of $Y$;
see \cite[X, Cor.\ 3.3]{SGA1}.
\end{rem}

% \begin{lem} \label{lem.sc-1} Suppose $Y$ is an iRreducible smooth simply
% connected affine variety and $U$ is an open subset of $Y$. If $Y \setminus U$  
% has codimension $\ge 2$ in $Y$ then $U$ is simply connected.
% \end{lem}
% 
% \begin{proof}
% We want to show that every finite \'etale cover $f \colon E \lra U$ 
% is split.  By our assumption on $Y$ it suffices to prove
% that $f$ can be extended to a finite \'etale cover over all of $Y$.
% 
% Since $U$ is smooth, so is $E$. Hence, $E$ is the integral
% closure of $U$ in $k(E)$, meaning that for every
% affine open $U_0 \subset U$, $\mathcal O_E(p^{-1}(U_0))$
% is the integral closure of
% $\mathcal O_{U}(U_0) \subset k(U)$ in $k(E)$. Let $R$ be
% the normal closure of the coordinate ring $k[Y]$ in $k(E)$,
% and put $Q = \Spec\,R$. By construction we have the following
% diagram:
% \[ \xymatrix{ E \;  \ar@{^{(}->}[r]  \ar@{->}[d] & 
% Q \ar@{->}[d]^{\pi} \cr
% U  \; \ar@{^{(}->}[r] & Y \, , } \]
% where $\pi$ is a finite map.  As $Y \setminus U$ 
% has codimension $\ge 2$, $\pi$ is \'etale in codimension two, 
% and therefore \'etale by the Zariski-Nagata
% purity theorem~\cite{nagata, zariski}. Thus $Q$ is
% the desired extension of $E$ to a finite \'etale cover over $Y$.
% \end{proof}

\begin{lem} \label{lem3.2} Let $G$ be a finite subgroup of $\GL(V)$
$\pi \colon V \to X = V \catq G$ be the categorical quotient 
map for the natural linear $G$-action on $V$. Then every morphism
$\alpha \colon V \lra V$ which commutes with $\pi$
is a translation by some $g \in G$.
\end{lem}

\begin{proof} 
Let $X_0$ be the principal stratum in $X$ and $U = \pi^{-1}(X_0)$.
Then $\pi_{|U} \colon U \lra X_0$ is a principal $G$-bundle;
in particular, there is a morphism $f \colon U \times_{X_0} U \lra G$ 
such that $u_2 = f(u_1, u_2) \cdot u_1$. Setting $g(u) = f(u, \alpha(u))$,
we obtain a morphism $g \colon U \lra G$ such that $\alpha(u) = g(u) \cdot u$.
Since $U$ is irreducible and $G$ is a finite group, 
the morphism $g$ is constant and the lemma follows.
\end{proof}

 \begin{lem} \label{lem2.2}
 Suppose a non-trivial finite subgroup $G \subset \GL(V)$ contains no 
 pseudo-reflections. Then $X = V \catq G$ is singular at 
 every point $x$ which does not lie in the principal Luna stratum.
 \end{lem}
 
  \begin{proof} 
  Let $\pi \colon V \lra X$ be the categorical quotient map, 
 $X_0$ be the principal stratum in $X$,
and $V_0 = \pi^{-1}(X_0) \subset V$. We want to show that
$\pi(v)$ is singular in $X$ for any $v \in V \subset V_0$. 
By the Luna Slice Theorem~\cite{luna} (see also~\cite[Section 6]{pv} 
or subsection~\ref{sect.slice} below)
it suffices to check that $V \catq \Stab(v)$ is singular 
at the origin. Since $v \not \in V_0$, $\Stab(v)$ is non-trivial.
The desired conclusion now follows from 
the Chevalley-Shephard-Todd theorem (cf., e.g., \cite[Theorem 8.1]{pv}), 
since we are assuming that
$G$ contains no pseudo-reflections (and hence, neither does $\Stab(v)$).
\end{proof}

\subsection{Proof of Proposition~\ref{prop3.2}}

\smallskip
(a) Let $X_0$ be the principal stratum in $X$
and $V_0 = \pi^{-1}(X_0) \subset V$. By Lemma~\ref{lem2.2}
$X_0$ is the smooth locus of $X$. Hence, 
$\sigma$ restricts to an automorphism of $X_0$. 

Now observe that the complement $V \setminus V_0$
is the union of $V^{g} = \{ v \in V \, | \, g(v) = v \}$,
as $g$ ranges over the non-identity
elements of $G$. Since no $g$ acts on $V$ as a pseudo-reflection,
each $V^g$ (with $g \ne e$) is a linear subspace of $V$ of
codimension $\ge 2$. Thus $V \setminus V_0$ has codimension $\ge 2$ in $V$.
Since $V$ is simply connected, Remark~\ref{rem.sc-1} tells us that
$V_0$ is simply connected as well.

Next we recall that $\pi_{|V_0} \colon V_0 \to X_0$ is 
a principal $G$-bundle
(in particular, it is \'etale).  By Lemma~\ref{lem.sc0}, applied 
to the map $f = \sigma \pi$ below
\[ \xymatrix{ V_0 \ar@{->}[d]^{\pi} \ar@{->}[r]^{\tilde{\sigma}} & 
V_0 \ar@{->}[d]^{\pi} \cr
X_0 \ar@{->}[r]_{\sigma} &  X_0  \, ,} \]
we see that $\sigma_{|X_0} \colon X_0 \to X_0$ 
lifts to a map $\tilde{\sigma} \colon V_0 \to V_0$. In other words,
\[ \sigma \pi (v) = \pi \tilde \sigma (v) \]
for every $v \in V_0$.  Since $V_0$ has codimension $\ge 2$ in $V$,
$\tilde \sigma$ extends to a morphism $V \lra V$, which we
will continue to denote by $\tilde \sigma$.
By Zariski density, the above equality holds for every $v \in V$.
In other words, $\tilde{\sigma}$ lifts $\sigma$ to $V$.

Similarly we can construct $\widetilde{\sigma^{-1}} \colon V \lra V$
lifting $\sigma^{-1}$ to $V$. Note that we are free to modify
$\tilde{\sigma}$ by composing it with a translation
$h \colon V \lra V$ for some $h \in G$ and that
the same is true of $\widetilde{\sigma^{-1}}$. Now observe that
$\alpha = \widetilde{\sigma^{-1}} \tilde \sigma$ satisfies $\pi =
\pi\alpha$. Thus by Lemma~\ref{lem3.2}
$\alpha$ is the translation morphism $g \colon V \lra V$ for some $g \in G$. After composing
$\tilde \sigma$ with $g^{-1}$ we may assume
$\tilde \sigma \widetilde{\sigma^{-1}} = \id_{V}$. In particular,
this shows that $\tilde{\sigma}$ is an automorphism of $V$.

\smallskip
(b) For any $g \in G$,
$\tilde \sigma g \tilde{\sigma}^{-1} \colon V \lra V$ commutes
with the quotient map $\pi$. Hence,
$\tau(g) := \tilde \sigma g \tilde{\sigma}^{-1} = h$ for some $h \in G$.
Clearly, $\tau$ is an automorphism of $G$, and part (b) follows.
\qed

\begin{remark} \label{rem.normalizer}
The statement of Proposition~\ref{prop3.2} does not assert that
the automorphism $\tilde \sigma \colon V \lra V$ is
linear. On the other hand, the automorphism $\tau$ of $G$
may actually be realized as conjugation by an element $a \in \GL(V)$
(normalizing $G$).  Indeed, differentiating
the equation $\tilde \sigma(gv) = \tau(g)\tilde
\sigma(v)$ at $v = 0$ we have $d\tilde \sigma_0 g = \tau(g) d\tilde
\sigma_0$. As $\tilde \sigma$ is an isomorphism, $a = d\tilde \sigma_0 \in
\GL(V)$ is the desired element.
            
Thus if a finite subgroup $G$ of $\GL(V)$ 
has no psedo-reflections and $N_{\GL(V)}(G)$ is 
generated by $G$ and $C_{\GL(V)}(G)$, then every Luna 
stratum in $V \catq G$ is intrinsic. In particular, this always
happens if every automorphism of $G$ is inner.
\end{remark}   

\begin{example} \label{ex.regular} Let $G$ be a finite group and $V = k[G]$
be its group algebra. Suppose $G$ has an automorphism $\tau$
and a subgroup $H$ such that $\tau(H)$ is not conjugate to $H$.
Then the Luna stratification in $X = V \catq G$ is intrinsic
but the stratum $X^{(H)}$ is not.
\end{example}

\begin{proof}
The existence of $\tau$ implies, in particular, that $|G| \ge 3$;
it is now easy to see that no $g \in G$ acts on $V$ as a pseudo-reflection.
Hence, the Luna stratification in $V \catq G$ is intrinsic. 

To see that $X^{(H)}$ is not intrinsic, let
$V^{\tau}$ be the vector space $V$, with the ``twisted" action
$g\star w := \tau(g)w$. As $V$ is the regular representation, there is a
$G$-equivariant linear isomorphism $\lambda \colon V \to V^{\tau}$.
The subgroup $H$ appears as an isotropy subgroup in $k[G/H]$,
which is a subrepresentation of $V$; hence, $H$ appears as an isotropy
subgroup in $V$. Clearly
$X = V \catq G \simeq V^{\tau} \catq G$ canonically, since the
polynomial invariants in $V^{\tau}$ are the same as in $V$.
Thus $\lambda$ descends to an automorphism $\sigma \colon X \lra X$.
By our construction $\sigma$ restricts to an isomorphism between
$X^{(H)}$ and $X^{(\tau(H))}$.
\end{proof}

It is easy to give examples of $G$, $\tau$ and $H$ satisfying 
the conditions of Example~\ref{ex.regular}. For instance,
we can take $G = (\mathbb{Z}/p \mathbb{Z})^r$ and $\tau$ 
to be a non-scalar element of $\GL_r(\mathbb{Z}/p \mathbb{Z})$.
Another family of examples can be constructed as follows.
Start with a non-trivial group $H$, set $G = H \times H$
and define an automorphism $\tau$ of $G$ by
$\tau(h_1, h_2) = (h_2, h_1)$.

\begin{remark} \label{rem.cs} Unfortunately, our proof of 
Theoprem~\ref{thm1} does not extend to representations 
of (infinite) reductive groups $G$. The main problem is 
that Remark~\ref{rem.sc-1}, which relies on the Zariski-Nagata
purity theorem, fails in this setting. 
If $U$ is an open subset of the affine space $\bbA^n$ 
whose complement has codimension $\ge 2$ then it is 
no longer true, in general, that every principal $G$-bundle can 
be extended from $U$ to $\bbA^n$, even if $n = 3$ 
and $G = \GL_n$; cf.~\cite[p. 124]{cs}. For this reason 
we cannot reproduce Proposition~\ref{prop3.2}(a) in
the setting of reductive groups. Our proof of 
Theorem~\ref{thm2} uses a rather different strategy,
based on studying the singularities of the Luna strata
in $V^r \catq G$; see Section~\ref{sect.strategy}.
\end{remark}

\section{Actions of reductive group}

In this section we collect several well-known definition and
results about actions of reductive groups on affine varieties,
in preparation for the proof of Theorem~\ref{thm2}.

\subsection{The Luna Slice Theorem}
\label{sect.slice}
Let $G \lra \GL(V)$ be a linear representation 
of a reductive group $G$ and $v \in V$ be
a point with a closed $G$-orbit. Then by Matsushima's theorem,
the stabilizer $H = \Stab(v)$ is a reductive subgroup of $G$.
Consequently, the $H$-subrepresentation $T_v(G \cdot v)$ of the natural
representation of $H$ on the tangent space $T_v(V)$ has
an $H$-invariant complement. We shall refer to this $H$-representation
as \emph{the slice representation} and denote it by $\Slice(v, V)$.
The Luna Slice Theorem asserts that the horizontal maps in the natural diagram
\[ \xymatrix{
 \Slice(v, V) *_H G \ar@{->}[d] \ar@{->}[r] & V \ar@{->}[d]^{\pi} \cr
 \Slice(v, V) \catq H \ar@{->}[r] & V \catq G } \]
are \'etale over $v$ and $\pi(v)$, respectively. For details,
see~\cite{luna} or~\cite[Section 6]{pv}. As an easy corollary, 
we obtain the following:

\begin{cor} \label{cor.slice1} Let $G$ be a reductive group, $G \lra \GL(V)$
be a linear representation, and $v \in V$, $H = \Stab(v) \subset G$
be as above. Then

\smallskip
{\rm(a)} the $G$-representation on $V$ is stable if and only
if the $H$-representation on $\Slice(v, V)$ is stable.

\smallskip
{\rm(b)} the $G$ representation on $V$ is generically free if and only if the
if the $H$-representation on $\Slice(v, V)$ is generically free.
\qed
\end{cor}

\subsection{Stability}

\begin{defn} \label{def.stable}
Let $G$ be a reductive group and $V$ be an affine $G$-variety.
A point $v \in V$ is called

\smallskip
{\em stable} if its orbit $G \cdot v$ is closed in $V$ and

\smallskip
{\em properly stable} if $v$ is stable and $\Stab_G(v)$ is finite.

\smallskip
\noindent
We shall say that the representation $V$ is

\smallskip
{\em stable} if a point $v \in V$ in general position
is stable,

\smallskip
{\em properly stable} if a point $v \in V$ in general position
is properly stable,

\smallskip
{\em generically free} if a point $v \in V$ in general position
has trivial stabilizer.
\end{defn}

Note that ``generically free" is not the same thing as ``having
trivial principal stabilizer". The reason is that when we talk
about the principal stabilizer, we are only interested in
$\Stab(v)$, where $v$ is a stable point. For example,
the natural action of the multiplicative group $\mathbb{G}_m$
on $V = \AA^1$ is generically free, but the principal isotropy
is all of $\mathbb{G}_m$, because the only stable point in $\AA^1$
is the origin. The precise relationship between these notions is
spelled out in the following lemma.

\begin{lem} \label{lem.stability}
Let $V$ be a linear representation of a reductive group $G$
and $\pi \colon V \to V \catq G$ be the categorical
quotient map.

\smallskip
{\rm(a)} If $\pi(v) \in (V \catq G)^{(\{ e \})}$ then $v$ is a properly
stable point in $V$.

\smallskip
{\rm(b)} The following conditions are equivalent.

\begin{itemize}
\item
 $V$ has trivial principal stabilizer,

\item
$V$ is generically free and properly stable,

\item
$V$ is generically free and stable.
\end{itemize}
\end{lem}

\begin{proof}
(a) Let $x = \pi(v) \in (V \catq G)^{(\{ e \})}$.
%  Assume $V$ has trivial principal stabilizer.
% Let $\pi \colon V \to V \catq G$
% be the categorical quotient map and $U \subset V \catq G$
% be the principal stratum.
We claim that $\pi^{-1}(x)$ is
a single $G$-orbit in $V$. Indeed, let $C = G \cdot v_0$ be the unique
closed orbit in $\pi^{-1}(x)$. Then $C$ is contained in the
closure of every orbit in $\pi^{-1}(x)$. On the other hand,
by our assumption $\Stab(v_0) = \{ e \}$;
hence, $\dim(C) = \dim(G)$, and $C$ cannot be contained
in the closure of any other $G$-orbit. This shows that $C = \pi^{-1}(x)$.
Thus every point in $\pi^{-1}(x)$ is stable (and hence, properly stable).
This proves part (a).  Part (b) is an immediate consequence of part (a).
\end{proof}

\subsection{The Hilbert-Mumford criterion}

Consider a linear $\GG_m$-representation on a vector space $V$.
Any such representation
can be diagonalized. That is, there is a basis $e_1, \dots, e_n$
of $V$, so that
$t \in \GG_m$ acts on $V$ by $t \cdot e_i \mapsto t^{d_i} e_i$.

In the sequel we shall use the following variant
of the Hilbert-Mumford criterion; see~\cite[Section 2.1]{mf}.

 \begin{thm} \label{thm.h-m} 
Consider a linear representation of a reductive group $G$ on
a vector space $V$. Then

\smallskip
{\rm(a)} $v \in V$ is properly stable for the action of
 $G$ if and only if it is properly stable for the action of
 every 1-dimensional subtorus $\GG_m \hookrightarrow G$.

\smallskip
{\rm(b)} In the above notations, $v = c_1 e_1 + \dots + c_n e_n \in V$
is properly stable for the action of $\GG_m$ if and only if there exist
$i, j  \in \{ 1, \dots, n \}$ such that $d_i < 0$, $d_j > 0$ and
$c_i, c_j \ne 0$.
\qed
 \end{thm}

We give several simple applications of this theorem below.

\begin{cor} \label{cor.h-m} Suppose $G$ is a reductive group,
$G \to \GL(V)$ is a linear representation,
$\pi_G \colon V \to V \catq G$ is the categorical quotient map,
$v \in V$ and, $\pi_G(v) \in (V \catq G)^{(\{ e \})}$.

\smallskip
{\rm(a)} If $H \subset G$ be a reductive subgroup and
$\pi_H \colon V \to V \catq G$ is the categorical quotient map for
the induced $H$-action on $V$ then $\pi_H(v) \in (V \catq H)^{(\{ e \})}$.

\smallskip
{\rm(b)} If $f \colon V' \lra V$ is a $G$-equivariant linear map,
$\pi' \colon V' \to V' \catq G$ is the categorical quotient
and $f(v') = v$ then $\pi'(v') \in (V' \catq G)^{(\{ e \})}$.
\end{cor}

\begin{proof}
Recall that, by definition, $\pi_G(v) \in V \catq G^{(\{ e \})}$ if
and only if into (i) $\Stab_G(v) = \{ e \}$ and
(ii) $v$ is properly stable for the $G$-action on $V$.

\smallskip
(a) We need to show that (i) and (ii) remain valid if $G$ is replaced by $H$.
In case of (i) this is obvious, and in case of (ii), this follows from
the Hilbert-Mumford criterion, since every 1-parameter subgroup of $H$
is also a 1-parameter subgroup of $G$.

\smallskip
(b) Again, we need to check that
$\Stab_{G}(v') = \{ e \}$ and
$v'$ is properly stable for the $G$-action on $V'$. The former
is obvious, and the latter follows from the Hilbert-Mumford criterion.
\end{proof}

\subsection{Reductive groups whose connected component is central}

Let $H$ be a reductive group whose connected component $H^0$
is central.  In particular, $H^0$ is abelian and hence, a torus.
Note that this class of groups includes both tori and all
finite groups (in the latter case $H^0 = \{ 1 \}$).

\begin{lem} \label{lem.central}
Let $H$ be a reductive group whose connected component $H^0$
is central and $\rho \colon H \to \GL(W)$
be a linear representation. Then

\smallskip
{\rm(a)} $\Stab(w) = \Ker(\rho)$ for $w \in W$ in general position.

\smallskip
{\rm(b)} Suppose $H$ has trivial principal stabilizer in $W^s$
for some $s \ge 1$.  Then $H$ has trivial principal
stabilizer in $W$.
\end{lem}

\begin{proof} (a) By~\cite[Theorem A]{richardson} the $H$-action
on $W$ has a stabilizer in general position.
That is, there is a subgroup $S \subset H$ and
an open subset $U \subset W$ such that $\Stab(u)$ is conjugate to $S$
for any $u \in U$. However, since $H^0$ is central in $H$, $S$
has only finitely many conjugates; denote them by $S = S_1, \dots, S_m$.
Then $U$ is contained in the union of finitely many linear subspaces
\[ U  \subset W^{S_1} \cup \dots \cup W^{S_m} \, . \]
Since $U$ is irreducible, we see that $U$ is contained in one of them,
say, $U \subset W^S$. Hence, $S \subset S_1, \dots, S_m$. Consequently,
$S$ is normal in $H$ and thus $\Stab(u) = S$ for every $u \in U$. 
This shows that $S = \Ker(\rho)$, and part (a) is proved.

\smallskip
(b) Let $\rho^s$ be the (diagonal) representation of $H$ on $W^s$.
Clearly $\Ker(\rho) = \Ker(\rho^s) = \{ e \}$. Part (a) now tells
us that since $\rho^s$ is generically free, so is $\rho$.

By Lemma~\ref{lem.stability} it now suffices to check that
the $H$-action $\rho$ on $W$ is properly stable. Let $\GG_m
\hookrightarrow H$ be a 1-dimensional subtorus. Diagonalize it
in the basis $e_1, \dots, e_n$ of $W$, so that it acts via
\[ t \colon e_i \lra t^{d_i} e_i \, . \]
Note that if we diagonalize the $\GG_m$-action on $W^s$ then
the same exponents $d_1, \dots, d_n$ will appear but each will
be repeated $s$ times. Thus by Theorem~\ref{thm.h-m},

\smallskip
\[ \begin{array}{ccc}
\text{$W^s$ is properly stable} & \Longleftrightarrow & 
\text{$d_i > 0$ and $d_j < 0$ for some $i, j \in \{ 1, \dots, n \}$} \\
 & & \Updownarrow \\
 & & \text{$W$ is properly stable,} 
\end{array} \]
and the lemma follows.
\end{proof}

\begin{remark} Note that both parts of Lemma~\ref{lem.central}
fail if we only assume that $H^0$ is a torus (but do not assume that
it is central in $H$). For example, both parts fail for
the natural action of the orthogonal group
$G = \OO_2(k) = \GG_m \rtimes \ZZ/2 \ZZ$ on $W = k^2$;
cf.~\cite[Example 2.5]{rv}.
\end{remark}

\subsection{Multiple representations}

 \begin{lem} \label{lem.slice2} Let $G$ be a reductive group and
$G \to \GL(W)$ be a linear representation.
Suppose that for some $r \ge 1$ an $r$-tuple $w = (w_1, \dots, w_r)$
is chosen so that $G \cdot w$ is closed in $W^r$ and that
$w_{d+1}, \dots, w_r$ are linear combinations of $w_1, \dots, w_d$ for
some $1 \le d \le r$.
Set $H = \Stab(w)$ and $v = (w_1, \dots, w_d) \in W^d$.  Then

 \smallskip
{\rm (a)} $\Stab_G (v) = H$.

 \smallskip
 {\rm(b)} $W^r$ has a $G$-subrepresentation $W_0$ such that $w \in W_0$
 and the natural projection $p \colon W^r \lra W^d$ onto the first $d$
 components restricts to an isomorphism between $W_0$ and $W^d$.

\smallskip
 {\rm(c)} $v$ has a closed orbit in $W^d$.

 \smallskip
 {\rm(d)} $\Slice(w, W^r) \simeq \Slice(v, W^d) \oplus W^{r-d}$,
 where $\simeq$ denotes equivalence of $H$-representations.
 \end{lem}

 \begin{proof} (a) is obvious. To prove (b),
 suppose $w_j = \sum_{i = 1}^d \alpha_{ij} w_j$
 for $j = d+1, \dots, r$. Then
 \[ \text{$W_0 = \{ (z_1, \dots, z_r) \, | \, z_j =
 \sum_{i = 1}^d \alpha_{ij} z_j$ for $j = d+1, \dots, r \; \}$.} \]
 has the desired properties. (c) follows from (b), since $G \cdot v$
 is the image of $G \cdot w$ under $p$.

\smallskip
 (d) Since $H$ is reductive,
 $W_0$ has an $H$-invariant complement $W_1$ in $W^r$, so
 that $W^r = W_0 \oplus W_1$. Since $W_0 \simeq W^d$, we conclude that
 $W_1 \simeq W^{r-d}$ (as an $H$-representation). The desired
 conclusion now follows from the fact that $p$ is an isomorphism between
 $W_0$ and $W^d$.
 \end{proof}

\begin{cor} \label{cor.slice3} Let $G$ be a reductive group and
$G \to \GL(W)$ be a linear representation of dimension $n$.
Then

\smallskip
{\rm(a)} The following are equivalent:
{\rm(i)} $W^r$ is stable for some $r \ge n$,
{\rm(ii)} $W^n$ is stable, and
{\rm(iii)} $W^s$ is stable for every $s \ge n$.

\smallskip
{\rm(b)} The following are equivalent:
{\rm(i)} $W^r$ is generically free for some $r \ge n$,
{\rm(ii)} $W^n$ is generically free, and
{\rm(iii)} $W^s$ is generically free for every $s \ge n$.

\smallskip
{\rm(c)} The following are equivalent:
{\rm(i)} $W^r$ has trivial principal stabilizer for some $r \ge n$,
{\rm(ii)} $W^n$ has trivial principal stabilizer, and
{\rm(iii)} $W^s$ has trivial principal stabilizer for
every $s \ge n$.
\end{cor}

\begin{proof} Suppose $r \ge n$ . Then for
$w = (w_1, \dots, w_r) \in W^r$ in general position,
$w_1, \dots, w_n$ span $W$. Keeping this in mind, we
see that

\smallskip
(a) the implication (i) $\Rightarrow$ (ii) follows from
Lemma~\ref{lem.slice2}(c) (with $d = n$)
and the implication (ii) $\Rightarrow$ (iii)
follows from Lemma~\ref{lem.slice2}(b) (again, with $d = n$).
(iii) $\Rightarrow$ (i) is obvious.

\smallskip
(b) follows from Lemma~\ref{lem.slice2}(a).
(c) follows from (a) and (b); see Lemma~\ref{lem.stability}(b).
\end{proof}

\section{Proof of Theorem~\ref{thm2}: the overall strategy}
\label{sect.strategy}

Let $G$ be a reductive algebraic group acting on a smooth
affine variety $V$ and $\pi \colon V \lra X = V \catq G$
be the categorical quotient map for this action.
We will say that $H$ is a {\em stabilizer subgroup} for
$V$ (or simply a {\em stabilizer subgroup}, if the reference
to $V$ is clear from the context), if $H = \Stab(v)$ for some
{\em stable} point $v \in V$. Clearly $H$ is a stabilizer subgroup
if and only if the Luna stratum $X^{(H)}$ associated
to its conjugacy class is non-empty.

\subsection{Strata that are singular along the boundary}

We will say that a Luna stratum $S$ is singular along its boundary
if the singular locus of its closure $\overline{S}$
is precisely $\overline{S} \setminus S$. The following lemma
is the starting point for our proof of Theorem~\ref{thm2}.

\begin{lem} \label{lem2.1}
Suppose every Luna stratum $S$ in $X = V \catq G$ is singular along
its boundary.  Then the Luna stratification in $X$ is intrinsic.
\end{lem}

\begin{proof}
Suppose $\sigma \colon X \lra X$ is an automorphism of $X$.
We want to show that for every Luna stratum $S \subset X$,
$\sigma$ restricts to an isomorphism between $S$ and another
Luna stratum $S'$. We will prove this by descending
induction on $\dim(S)$.

The base case is the principal stratum $X_0$. By our assumption
this stratum is precisely the smooth locus of $X_0$; thus
$\sigma(X_0) = X_0$.

For the induction step, suppose $S \ne X_0$. We want to show that
$\sigma(S)$ is again a Luna stratum in $X$.
As we remarked in the Introduction, the (finitely many) Luna strata
are partially ordered; $S \preceq T$ if $S$ lies in the closure of $T$.
Let $T$ be a minimal
stratum with the property that $S \prec T$.  Then $\dim(T) > \dim(S)$
and thus, by the induction assumption,
$\sigma$ restricts to an isomorphism between
$T$ and another stratum, $T'$ and hence, also to an isomorphism
between their closures. Denote the maximal strata contained in
$\overline{T} \setminus T$ by $S_1 = S$, $S_2, \dots, S_n$ and
the maximal strata contained in
$\overline{T'} \setminus T'$ by $S_1'$, $S_2', \dots, S_n'$.
Then $\overline{S_1}, \dots, \overline{S_n}$
are the irreducible components of the singular locus
$\overline{T} \setminus T$ of $\overline{T}$ and
similarly $\overline{S_1'}, \dots, \overline{S_n'}$
are the irreducible components of the singular locus
$\overline{T'} \setminus T'$ of $\overline{T'}$.
Thus $\sigma$ takes each $\overline{S_i}$ isomorphically to some
$\overline{S_j}'$.  In
particular, $\sigma$ restricts to an isomorphism between
$\overline{S}$ and $\overline{S_j}'$ for some
$j = 1, \dots, s$. Finally, restricting $\sigma$ to
an isomorphism between the smooth loci of
$\overline{S}$ and $\overline{S_j}'$, we see that
$\sigma$ takes $S$ to the stratum $S_j'$, as desired.
\end{proof}

\subsection{Acceptable families of representations}

\begin{defn} \label{def.acceptable}
We shall call a family $\Lambda$ of finite-dimensional linear 
representations $G \to \GL(V)$ of reductive (but not necessarliy
connected) algebraic groups {\em acceptable} if it satisfies 
the following two conditions.

\smallskip
(i) If $G \lra \GL(V)$ is in $\Lambda$ then for every stabilizer subgroup
$H$ in $G$, the induced representation $N_G(H) \to \GL(V^H)$ is again in
$\Lambda$, and

\smallskip
(ii) For every representation $G \lra \GL(V)$ in $\Lambda$,
the principal stratum in $V \catq G$ is singular along its boundary.
\end{defn}

\begin{prop} \label{prop.strategy}
Suppose a linear representation $G \lra \GL(V)$ belongs
to an acceptable family. Then every Luna stratum in $X = V \catq G$
is singular along its boundary. In particular, the Luna stratification
in $X$ is intrinsic.
\end{prop}

\begin{proof} The second assertion follows from the first 
by Lemma~\ref{lem2.1}. Hence, we only need to show that every
Luna stratum in $X$ is singular along its boundary. 
Let $\pi \colon V \to V \catq G = X$ be the categorical
quotient map and let $S = X^{(H)}$ be a Luna stratum.
Choose $p \in \overline{S} \setminus S$, say, $p \in X^{(K)}$, where
$K$ is a (reductive) stabilizer subgroup and $H \subsetneq K$.
Our goal is to show that $\overline{S}$ is singular at $p$.
We will argue by contradiction.  Assume, to the contrary,
that $\overline{S}$ is smooth at $p$.

Let $N = N_G(H)$ be the normalizer
of $H$ in $G$ and write the surjective
map $\pi_{| \, V^H} \colon V^H \lra \overline{S}$ as a composition
\[ \xymatrix{ V^H \ar@{->}[d]^{\pi_N} V^H \ar@
/^2pc/[dd]^{\pi} \cr
  V^H \catq N \ar@{->}[d]^{n} \cr
  \overline{S},} \]
where $\pi_N$ is the categorical quotient map for the $N$-action
on $V^H$ and $n$ is the normalization map for $\overline{S}$.

Let $v \in V^H$ be an $N$-stable point with stabilizer $K$
such that $\pi(v) = p$ and let $q = \pi_N(v)$. (Note that by Luna's
criterion, $v$ is $N$-stable if and only if it is $G$-stable;
see~\cite[Theorem 6.17]{pv}.) Recall that we are assuming
that the $G$-representation on $V$ belongs to an acceptable family
$\Lambda$.  Consequently, the $N$-representati
on on $V^H$
also belongs to $\Lambda$, and thus the smooth locus of $V^H \catq N$
is precisely the principal statum for the $N$-action on $V^H$.
In other words, if $q$ does not lie in the principal stratum
in $V^H \catq N$ then $q$ is a singular point of
$V^H \catq N$. Since $n$ is the normalization map and $n(q) = p$,
this implies that $p$ is a singular point of $\overline{S}$,
a contradiction.

We may thus assume that $q$ lies in the principal stratum $U$
of $V^H \catq N$.  Since we are assuming that
$p$ is a smooth point of $\overline{S}$, the normalization map
$n$ is an isomorphism between Zariski open neighborhoods 
of $p$ and $q$.  Since $\pi_N^{-1}(U) \to U$ is a
principal $N$-bundle, it follows that the differential
$d \pi_v$ maps  $T_{v}(V^H)$ surjectively onto $T_{p}(\overline{S})$.

We will now show that this is impossible.
Indeed, since the quotient map
$\pi \colon V \lra X$ is $K$-equivariant
(where $K$ acts trivially on $X$) and
$v$ is fixed by $K$, the differential
$d \pi_v \colon T_{v}(V) \lra T_p(X)$ is
a $K$-equivariant linear map (where $K$ acts trivially
on $T_p(X)$). Consequently, $d \pi_v$ sends every
non-trivial irreducible $K$-subrepresentation
of $T_{v}(V)$ to $0$.
Since $V$ is smooth, we have $(T_{v}(V))^K
= T_{v}(V^K)$ and therefore $d\pi_{v}$ maps $T_{v}(V^K)$ onto
$T_p(\overline{S})$. On the other hand, since
$\pi(V^K) = \overline{X^{(K)}}$, we conclude that
\begin{equation} \label{e2.1}
T_p(\overline{X^{(K)}}) \supset T_p(\overline{S}) \, .
\end{equation}
Now recall that we are assuming that $p$ is a smooth point of $\overline{S}$.
Moreover, since $p \in X^{(K)}$, it is also a smooth point
of $\overline{X^{(K)}}$. Thus~\eqref{e2.1}
implies $\dim X^{(K)} \ge \dim S$, contradicting the fact that
$X^{(K)}$ lies in $\overline{S} \setminus S$.
\end{proof}

We now record a corollary of the above argument for future reference.

\begin{cor} \label{cor1.acceptable} Suppose $G$ is a reductive group
and $G \lra \GL(V)$ is a linear representation with principal
isotropy $H$. Let $N = N_G(H)$. Then the natural map
\[ \xymatrix{V \catq N \ar@{->}[d]^{n} \cr
  V \catq G } \]
is an isomorphism, which identifies the principal Luna stratum
in $V^H \catq N$ with the principal Luna stratum in $V \catq G$.
\end{cor}

\begin{proof}
The fact that $n$ is an isomorphism is proved in~\cite[Corollary 4.4]{lr}.

To show that $n$ identifies the principal strata in $V \catq N$ and
$V \catq G$, let $\pi \colon V \to V \catq G$ and 
$\pi_N \colon V \to V \catq N$
be the categorical quotient maps.
Choose $p \in V \catq N$ and set $q = n(p)$. Let
$v \in V$ be a point in the (unique) closed $N$-orbit in $\pi_N^{-1}(p)$.
By Luna's criterion, the $G$-orbit of $v$ is also closed;
cf. see~\cite[Theorem 6.17]{pv}.  Our goal is to show that
\[ \Stab_G(v) = H \Longleftrightarrow \Stab_N(v) = H \, . \]
The $\Rightarrow$ direction is obvious, so suppose $\Stab_N(v) = H$,
(i.e., $q$ lies in the principal stratum in $V^H \catq N$) and
$\Stab_G(v) = K$. We want to show that $K = H$. Assume the contrary:
$H \subsetneq K$. Since $p$ lies in the principal stratum in $V \catq N$,
$d \pi_N$ maps $T_v(V)$ surjectively to $T_p(V \catq N)$. Since
$n$ is an isomorphism and $\pi = n \circ \pi_N$, we see that
$d \pi_v$ maps $T_v(V)$ surjectively onto $T_v(V \catq G)$.
On the other hand, in the proof of Proposition~\ref{prop.strategy}
we showed that this is impossible if $H \subsetneq K$.
\end{proof}

\subsection{The strategy}
\label{rem.strategy}

Our proof of Theorem~\ref{thm2} will be based on showing that each
of the families of representations in parts (a), (b) and (c)
is acceptable, i.e., satisfies conditions (i) and (ii) of
Definition~\ref{def.acceptable}; the desired conclusion will
then follow from Lemma~\ref{lem2.1} and Proposition~\ref{prop.strategy}.
In particular, our argument will also show that for
the representations $V$ considered
in Theorem~\ref{thm2}, each Luna stratum in $V \catq G$ is singular
along its boundary.

\subsection{A toy example}
To illustrate this strategy, we will apply it to
the following simple example. We will say that
a linear representation $G \lra \GL(V)$ has has {\em the
codimension $2$ property} if $\dim(V^A) - \dim(V^B) \ne 1$ for every
pair of subgroups $A \triangleleft B \le G$ (here
$A$ is normal in $B$). 

\begin{prop} Let $\Lambda$ be the family of representations
$\phi \colon G \lra \GL(V)$, where $G$ is finite and $\phi$ has
the codimension $2$ property. Then $\Lambda$ is acceptable.
\end{prop}

\begin{proof}
Condition (ii) of Definition~\ref{def.acceptable} follows from
Lemma~\ref{lem2.2} and condition (i) is immediate from the definition.
\end{proof}

We conclude that if
$G$ be a finite group and $G \lra \GL(V)$ be a representation with the
codimension $2$ property. Then

\smallskip
(a) every Luna stratum in $V \catq G$
is singular along its boundary, and 

\smallskip
(b) the Luna stratification in $V \catq G$ is intrinsic.  

\smallskip
Clearly no element of $G$ can act on $V$ as
a pseudo-reflection (to see this, set $B = \left< g \right>$
and $A = \{ 1_G \}$); thus (b) is a special case of
Theorem~\ref{thm1}. On the other hand, (a) does not follow
from either the statement or the proof of Theorem~\ref{thm1}
(at least, not from our proof).

The following example which shows that there are many 
representations with the codimension $2$ property.
Part (a) will be used again in the sequel.

\begin{example} \label{ex.orthogonal1} (a) Let $G \lra \GL(W)$ be a linear
representation, leaving invariant a nondegenerate bilinear form $b$
on $W$. Then for any reductive subgroup $H$ of $G$ the restriction of
$b$ to $W^H$ is again non-degenerate.

\smallskip
(b) Every symplectic representation of a finite
group has the codimension $2$ property.
\end{example}

Recall that a (not necessarily symmetric)
a bilinear form $b$ on $W$ is called {\em non-degenerate} 
if for every $0 \ne a \in W$ the linear form $l_a \colon w \mapsto b(a, w)$
is not identically $0$ on $W$. 

\begin{proof} (a) Suppose for some $a \in W^H$ the linear form
$l_a$ given by $w \mapsto b(a, w)$ is identically zero on $W^H$.
We want to show that $a = 0$. Since $b$ is non-degenerate, 
it suffices to show that $l_a$ is identically zero on all of $W$. 

Since $a \in W^H$, $l_a$ is a linear invariant for the $H$-action on $W$.  
By Schur's lemma, $l_a$ vanishes on every non-trivial irreducible
$H$-subrepresentation of $W$. On the other hand, we are assuming that
$l_a$ also
vanishes on $W^H$, which contains every trivial irreducible
$H$-subrepresentation 
of $W$. Thus $l_a$ is identically zero on $W$, as claimed.

\smallskip
(b) Suppose $G$ preserves a symplectic form $b$ on $W$.
%  \colon \bigwedge^2(W) \lra k$.
Then part (a) tells us that $W^H$ is even-dimensional for every reductive 
subgroup $H \subset G$. The codimension $2$ property is now immediate
from the definition.
\end{proof}

\section{Non-coregular representations}

The main difficulty in implementing the strategy
outlined in section~\ref{rem.strategy} is in checking condition
(ii) of Definition~\ref{def.acceptable}. That is, given a linear
representation $G \to \GL(V)$ of a reductive group $G$ on
a vector space $V$ and a stable point $v \in V$, we want
to show that $V \catq G$ is singular at $\pi(v)$.
The Luna Slice Theorem reduces this problem to checking that
$\Slice(v, V) \catq H$ is singular at $\pi_H(0)$, where
$H = \Stab(v)$, and $\pi_H \colon \Slice(v, V) \to \Slice(v, V) \catq H$
is the categorical quotient map. In other words, we are reduced
to a problem of the same type, where $v = 0$ and $G = H$.
% (Note however, that if the original representation
% of $G$ on $V$ belongs to an acceptable family $\Lambda$,
% the slice representation of $H$ on $\Slice(v, V)$ may
% not belong to $\Lambda$. We shall not be concerned with this
% complication for the time being.)

Linear representation $G \to \GL(V)$, with the property that
$V \catq G$ is smooth at $\pi(0)$ are called {\em coregular}.
Here $G$ is assumed to be reductive and $\pi$ is the categorical
quotient map $V \to V \catq G$. It is easy to see 
(cf., e.g.,~\cite[Proposition 4.11]{pv}) that $V \catq G$ 
is smooth at $\pi(0)$ if and only if it is smooth everywhere
if and only if it is an affine space $\bbA^d$ for some $d \ge 1$.
Coregular representations have been extensively studied;
for a survey of this topic and further references,
see~\cite[Section 8]{pv}. Thus in order to implement the strategy
for proving Theorem~\ref{thm2} outlined in section~\ref{rem.strategy}
we need a large family of representations that are known
not to be coregular.  The purpose of this section is to prove
Proposition~\ref{prop5.1} below, which exhibits such a family.

\begin{prop}\label{prop5.1}
Let $G$ be a reductive group and $V_1, V_2$ be linear representations of
$G$, such that
$V_1$ has trivial principal stabilizer and $V_2$ is not
fixed pointed (see Definition~\ref{def.fixed-pointed} below).
Then $V_1 \times V_2$ is not coregular. That is, 
$(V_1 \times V_2) \catq G$ is singular.
\end{prop}

We begin with preliminary results about fixed pointed
representations, which will be used in the proof and
in subsequent applications of Proposition~\ref{prop5.1}.

\subsection{Fixed pointed representations}

\begin{defn} \label{def.fixed-pointed} Following Bass
and Haboush~\cite{bh}, we will say that
a linear representation $G \to \GL(V)$ of a reductive group $G$
is {\em fixed pointed} if the natural $G$-equivariant
projection $\pi \colon V \lra V^G$ is
the categorical quotient for the $G$-action on $V$.
(Note that the projection $\pi$ is sometimes called
the Reynolds operator; cf., e.g., \cite[Section 3.4]{pv}.)
\end{defn}

\begin{lem} \label{lem.fixed-pointed1}
Let $G$ be a reductive group and $\rho \colon G \to \GL(V)$ be a
linear representation. Then the following conditions are equivalent.

\smallskip
{\rm(a)} $\rho$ is fixed pointed,

\smallskip
{\rm(b)} The null cone $\NC(V)$ is {\rm(}scheme-theoretically{\rm)} a vector space,

\smallskip
{\rm(c)} The null cone $\NC(V)$ is {\rm(}scheme-theoretically{\rm)} smooth.

\smallskip
\noindent
If {\rm(a)}, {\rm(b)} and {\rm(c)} hold then $V = NC(V) \oplus V^G$.
\end{lem}

\begin{proof}
(a) $\Rightarrow$ (b).  Let $\pi \colon V \lra V \catq G$ be
the categorical quotient map. If $\rho$ is fixed pointed then
$\pi$ is a linear projection, so clearly $\NC(V) = \pi^{-1}(0)$
and $V = NC(V) \oplus V^G$.

(b) $\Rightarrow$ (c) is obvious, since a vector space is smooth.

(c) $\Rightarrow$ (a). Assume that $\NC(V)$ is scheme-theoretically
smooth.  Since $G$ is reductive, we may write
$V = V^G \oplus W$ for some $G$-invariant subspace $W$.
Then clearly $W^G = (0)$, $\NC(W) = \NC(V)$ and
$V \catq G = V^G \times W \catq G$, where the quotient map
$\pi_V$ sends $(v, w) \in V = V^G \oplus W$ to 
$(v, \pi_W(w)) \in V^G \times W \catq G$.

Thus, after replacing $V$ by $W$, we may assume without loss
of generality that $V^G = (0)$.  Our goal is then to show that
$V$ is fixed pointed, which, in this case, means that
$V \catq G$ is a single point, or equivalently,
\[ \NC(V) = V \, . \]
Indeed, assume
the contrary. Then $\NC(V)$ is cut out by the equations $f = 0$,
as $f$ ranges over the homogeneous elements of $k[V]^G$.

Note that no nonzero homogeneous element $f \in k[V]^G$ can
be of degree $1$. Indeed, if there were
a non-zero $G$-invariant linear function $f \colon V \lra k$
then $V$ would contain a copy of the trivial representation,
contradicting $V^G = (0)$.
We thus conclude that $\NC(V)$ is a subscheme of $V$ cut out
by a (possibly empty) collection of a homogeneous polynomials 
of degree $\ge 2$. In particular, the tangent space to $\NC(V)$
at $0$ coincides with all of $T_0(V)$. Since we are assuming that
$\NC(V)$ is (scheme-theoretically) smooth, this is only possible if
$\NC(V) = V$.
\end{proof}

In the sequel we will primarily be interested in representations that
are {\em not} fixed pointed. Two examples are given below.

\begin{example} \label{ex.fixed-pointed2}
No nontrivial stable
representation $G \to \GL(V)$ can be fixed pointed.
Indeed, in a fixed-pointed representation, the only
stable points are those in $V^G$.
\end{example}

\begin{example} \label{ex.orthogonal2}
A non-trivial orthogonal representation
$H \lra \OO_{\phi}(L)$ of a reductive group $H$ on
a vector space $L$ (preserving a non-degenerate
quadratic form $\phi$) is not fixed pointed.
\end{example}

\begin{proof} Assume the contrary. Then
\begin{equation} \label{e.o2a}
L = \NC(L) \oplus L^H \, \end{equation}
where $\NC(L)$ is the null-cone of $L$.  Example~\ref{ex.orthogonal1}(a),
tells us that
$\phi$ restricts to a non-degenrate quadratic form on $L^H$. Hence,
\[ L = (L^H)^{\perp} \oplus L^H \, . \]
Now observe that $L^H$ has a unique $H$-invariant complement in $L$
(namely, the direct sum of all non-trivial $H$-subrepresentations in $L$).
We thus conclude that $\NC(L) = (L^H)^{\perp}$.
Since $\phi$ is non-degenerate on $L$ and $L^H$,
it is also non-degenerate on $(L^H)^{\perp}$. Note that
since the $H$-action on $L$ is non-trivial, $L^H \ne L$
and thus $\NC(L) = (L^H)^{\perp} \ne (0)$. In particular,
the $H$-invariant regular function $L \lra k$ given
by $x \lra \phi(x, x)$ is not constant on $(L^H)^{\perp}$.
On the other hand, every $H$-invariant regular function
on $\NC(L)$ has to be constant. This contradiction
shows that $\phi$ is not fixed pointed.
\end{proof}

\subsection{Proof of Proposition~\ref{prop5.1}}

Set $V = V_1 \times V_2$. Assume the contrary:
$V \catq G$ is smooth, i.e.. is isomorphic to an affine space $\AA^d$.
Let $\pi\colon V\to V\catq G$, $\pi_1 \colon V_1 \to V_1\catq G$,
and $\pi_2 \colon V_2 \to V_2 \catq G$ be the categorical
quotient maps. We will denote the projection $V \to V_i$ by $p_i$ and
the induced morphism $V \catq G \to V_i \catq G$ by $\overline{p_i}$.

Let $T = \bbG_m \times \bbG_m$ 
be a two-dimensional torus acting on $V = V_1 \times V_2$
by \[ (s, t) \colon (v_1, v_2) \to (s v_2 , t v_2) \, . \]
This action commutes with the $G$-action on $V$ and hence,
decends to $V \catq G$. Clearly the $T$-fixed point $\pi(0, 0)$
lies in the closure of every other $T$-orbit in $V \catq G \simeq \AA^d$.
Thus by \cite[Corollary 10.6]{bh}, the $T$-action on $V \catq G$
is isomorphic to a linear action. That is, we may assume that
$V \catq G$ is a vector space with a linear action of $T$.
This identifies $V_1 \catq G$ and $V_2 \catq G$ with the $T$-invariant
linear subspaces
$(V \catq G)^{\{ 1 \} \times \bbG_m}$ and 
$(V \catq G)^{\bbG_m \times \{ 1 \}}$ of $V \catq G$ respectively.
In particular, $V_i \catq G$ is smooth for $i = 1, 2$. Moreover,
\[ \overline{p_1} (\overline{v}) =
\lim_{t \to 0} (1, t) \cdot \overline{v} \]
and
\[ \overline{p_2} (\overline{v}) =
\lim_{s \to 0} (s, 1) \cdot \overline{v} \]
are $T$-equivariant linear projections $V \catq G \to V_1 \catq G$ 
and $V \catq G \to V_2 \catq G$ respectively and
\[ \overline{p} = (\overline{p_1}, \overline{p_2})
 \colon V \catq G \to (V_1 \catq G) \times (V_2 \catq G) \]
is a smooth map.

Consider the commutative diagram
\[   \xymatrix{ V \ar@{->}[d]^{\pi}  \ar@{=}[r] & V_1 \times V_2  
\ar@{->}[ddl]^{\pi_1 \times \pi_2} \cr
 V \catq G \ar@{->}[d]^{\overline{p}} &  
% V_1 \times V_2 \ar@{->}[dl]^{\pi_1 \times \pi_2} 
\cr
 V_1 \catq G \times V_2 \catq G   &   }
\]
Choose a stable point $v_1 \in V_1$ such that $\pi_1(v_1)$ lies in the
principal stratum of $V_1 \catq G$ ($V_1$ has a dense open subset consisting
of such points; cf. Lemma~\ref{lem.stability}) and
let $x = (\pi_1(v_1), \pi_2(0)) \in V_1 \catq G \times V_2 \catq G$.
The (scheme-theoretic) preimage of $x$ under
the map $\pi_1 \times \pi_2$ is clearly $G \cdot v_1 \times \NC(V_2)$,
where $\NC(V_2)$ is the null cone in $V_2$.
Since we are assuming that the $G$-action on $V_2$ is not fixed pointed,
$\NC(V_2)$ is singular; cf. Lemma~\ref{lem.fixed-pointed1}. Thus
\[ \text{ $(\pi_1 \times \pi_2)^{-1}(x)$ is singular.} \]
On the other hand, as we saw above,
the map $\overline{p}$ is smooth. The map $\pi$ is
smooth over the principal Luna stratum in $V \catq G$.
By Corollary~\ref{cor.h-m}(b), every point in $p^{-1}(x)$ lies in
the principal stratum of $V \catq G$.
Hence, the composition map $\pi_1 \times \pi_2 = \overline{p}
\pi \colon V \lra V_1 \catq G \times V_2 \catq G$ is smooth
over some Zariski open neighborhood of $x$. Consequently,
\[ \text{
$(\pi_1 \times \pi_2)^{-1}(x)$ is smooth.} \]
This contradiction shows that $V \catq G$ cannot be smooth,
thus completing the proof of Proposition~\ref{prop5.1}.
\qed

\begin{cor} \label{cor1-to-prop5.1}
Suppose $V_1$, $V_2$ and $V_3$ are three linear representations
of a reductive group $G$, where $V_1$ has trivial principal
stabilizer, $V_2$ is not fixed pointed, and $V_3$ is arbitrary. Then
$V_1 \oplus V_2 \oplus V_3$ is not coregular.
\end{cor}

\begin{proof} The $G$-representation $V_1' = V_1 \oplus V_3$
has trivial principal stabilizer; see Corollary~\ref{cor.h-m}(a).
Now apply Proposition~\ref{prop5.1} to $V_1' \oplus V_2$.
\end{proof}

\section{Proof of Theorem~\ref{thm2}(a) and (b)}
\label{sec:REDUCTIVE}

We will follow the strategy outlined in section~\ref{rem.strategy}
by exhibiting acceptable families $\Lambda_a$ and $\Lambda_b$
which include the representations in parts (a) and (b)
of Theorem~\ref{thm2}, respectively.

Elements of $\Lambda_a$ are representations
of the form $V = W^r$, where $r \ge 2 \dim(W)$
and $G \lra \GL(W)$ is a representation of a reductive group $G$.

Elements of $\Lambda_b$ are representations
of the form $V = W^r$, where $r \ge \dim(W) + 1$
and $G \lra \OO(W)$ is an orthogonal representation
of a reductive group $G$. That is, $G$ preserves some non-degenerate
quadratic form on $W$.

In view of Proposition~\ref{prop.strategy}
it suffices to show that $\Lambda_a$ and $\Lambda_b$ are
acceptable families. We begin by checking condition (i)
of Definition~\ref{def.acceptable}.

\smallskip
(a) For $\Lambda_a$ this is
quite easy. Suppose $V = W^r$ is in $\Lambda_a$. Then
$V^H = (W^H)^r$ is again in $\Lambda_a$, because
$r \ge 2 \dim(W) \ge 2 \dim(W^H)$.

\smallskip
(b) Suppose $V = W^r$, where $W$ is an orthogonal representation of $G$
and $r \ge \dim(W) + 1$. Once again, $V^H = (W^H)^r$,
where $r \ge \dim(W) + 1 \ge \dim(W^H) + 1$. Moreover,
in view of Example~\ref{ex.orthogonal1}(a),
the $N_G(H)$-representation on $W^H$
is orthogonal. Thus $V^H$ belongs to $\Lambda_b$, as claimed.

\smallskip
It remains to show that $\Lambda_a$ and $\Lambda_b$ satisfy
condition (ii) of Definition~\ref{def.acceptable}. That is, given
a representation $V = W^r$ in $\Lambda_a$ or $\Lambda_b$, we want
to show that $X = V \catq G$ is singular at every point $x$ away from
the principal stratum. We begin with two reductions.

\smallskip
First we claim that we may assume without loss of generality
that $V$ has trivial principal stabilizer.  Indeed, suppose
the principal stabilizer in $V$ is $H \subset G$.
Set $N = N_G(H)$ and $\overline{N} = N/H$.
Then by Corollary~\ref{cor1.acceptable}
\[ \xymatrix{ V^H \catq \overline{N} \ar@{=}[r] & V^H \catq N \ar@{->}[d]^{n} 
\cr & X} \]
is an isomorphism which takes the principal stratum
in $V^H \catq N$ to the principal stratum in $X$.
Thus it suffices to prove that $V^H \catq N$ is singular
away from its principal stratum. As we just showed, the representation
$N \to \GL(V^H)$ lies in $\Lambda_a$ in
part (a) and in $\Lambda_b$ in part (b). Hence, so does
$\overline{N} \to \GL(V^H)$. Since the latter representation
has trivial principal stabilizer, this proves the claim.

\smallskip
 From now on, we will assume that the principal stabilizer
subgroup of $G$ in $V = W^r$ (but not necessarily in $W$) is $\{ e \}$.
Suppose $x$ is represented by an element $w = (w_1, \dots, w_r) \in W^r$
whose $G$-orbit is closed. Let $H = \Stab(v)$. Note that since
$x$ does not lie in the principal stratum in $V \catq G$, $H \neq \{ e \}$.
After permuting the components of $W^r$ if necessary, we may assume
that $w_{n+1}, \dots, w_r$ are linear combinations of $w_1, \dots, w_n$.

\smallskip
Recall that by the Luna Slice Theorem $(W^r \catq G, x)$ is
\'etale isomorphic to \[ (\Slice(w, W^r) \catq H, \pi(0)) \, ; \]
cf. e.g.,~\cite[Section 6]{pv}. Thus it suffices to prove that
$\Slice(w, W^r) \catq H$ is singular.

Let $v = (w_1, \dots, w_n) \in W^n$. By Lemma~\ref{lem.slice2}
$\Stab(v) = H$, $G \cdot v$ is closed in $W^n$, and
\[ \Slice(w, W^r) \simeq \Slice(v, W^n) \oplus W^{r-n} \, . \]
Recall that we are assuming that the $G$-action (and hence, the $H$-action)
on $W^r$ has trivial principal stabilizer; cf. Corollary~\ref{cor.h-m}(a).
Since $r \ge n$ (both in part (a) and in part (b)),
Corollary~\ref{cor.slice3}(c) tells us that the $H$-action on $W^n$
also has trivial principal stabilizer. By Corollary~\ref{cor.slice1}
this implies that the $H$-action on $\Slice(v, W^n)$ 
also has trivial principal stabilizer.

We will now consider the families $\Lambda_a$ and $\Lambda_b$ separately.

\smallskip
(a) Suppose $V$ is in $\Lambda_a$. Recall that we are assuming
$r \ge 2n$, i.e. $r-n \ge n$. Thus
by Corollary~\ref{cor.slice3}(c) $W^{r-n}$
also has trivial principal stabilizer. Consequently, it is not
fixed pointed; cf. Example~\ref{ex.fixed-pointed2}.
By Proposition~\ref{prop5.1} the $H$-representation
\[ \Slice(w, W^r) \simeq
\underbrace{\Slice(v, W^n)}_{\text{\tiny trivial principal stabilizer}}
\oplus
\underbrace{W^{r-n}}_{\text{\tiny not fixed pointed}} \]
is not coregular.

\smallskip
(b) Since $G$ preserves the non-degenerate
quadratic form $q \oplus \dots \oplus q$ ($n-r$ times) on $W^{r-n}$
and $r \ge n+1$, Example~\ref{ex.orthogonal2} tells us that
the $H$-representation on $W^{r-n}$ is not fixed pointed.
Proposition~\ref{prop5.1} now tells us that the $H$-representation
\[ \Slice(w, W^r) \simeq
\underbrace{\Slice(v, W^n)}_{\text{\tiny trivial principal stabilizer}}
\oplus
\underbrace{W^{r-n}}_{\text{\tiny not fixed pointed}} \]
is not coregular.
\qed

\section{Proof of Theorem~\ref{thm2}(c)}

Once again, we will follow the strategy outlined
in section~\ref{rem.strategy}, by defining a suitable
family $\Lambda_c$ of representations of reductive groups
and then checking that $\Lambda_c$ is acceptable.
In the context of Theorem~\ref{thm2}(c) the natural
candidate for $\Lambda_c$ is the family of representations 
of the form $W^r$, where $W = \Lie(G)$ is the adjoint 
representation for some reductive group $G$.
Unfortunately, this family is not acceptable, because 
it does not satisfy condition (i) of 
Definition~\ref{def.acceptable}. To make the argument 
go through, we need to consider a slightly larger family
which we now proceed to define.

\begin{defn} \label{def.almostadjoint}
Let $G$ be a reductive group.  We will say that 
a linear representation $\rho \colon G \to \GL(W)$
is {\em almost adjoint} if $\Ker(\rho)$ contains 
a normal closed subgroup $K$ of $G$ such that $W$ 
is isomorphic to the Lie algebra of $G/K$
and $\rho$ can be written as a composition 
\begin{equation} \label{e.almostadjoint}
\rho \colon G \to G/K \stackrel{\Ad_{G/K}} \to \GL(W) \, , 
\end{equation}
where $G \to G/K$ is the natural quotient map
and $\Ad_{G/K}$ is the adjoint representation.
Note that the groups $G$, $K$ or $G/K$ are assumed 
to be reductive but not necessarily connected.
\end{defn}

We are now ready to define $\Lambda_c$. Fix an integer
$r \ge 3$ and let $\Lambda_c$ be 
the family of representations of the form
$W^r$, where $G \lra \GL(W)$ is almost adjoint.
Following the strategy of section~\ref{rem.strategy}, 
in order to prove Theorem~\ref{thm2}(c), it suffices
to check that $\Lambda_c$ is an acceptable family.

We begin by checking condition (i) of Definition~\ref{def.acceptable}.
Suppose $\rho \colon G \lra \GL(W)$ is an almost adjoint representation,
with $K \triangleleft \, G$ as in~\eqref{e.almostadjoint}
and $H \subset G$ is a stabilizer subgroup. Since $(W^r)^H = (W^H)^r$,
it suffices to show that the natural representation of
the normalizer $N = N_G(H)$ on $W^H$ is again almost adjoint.

Note that $H$ contains $K$; since $N/K = N_{G/K}(H/K)$, we may,
after replacing $G$ by $G/K$, assume without loss of 
generality that $K = \{ e \}$, i.e., $W = \Lie(G)$ is the Lie algebra 
of $G$ and $\rho$ is the adjoint representation.

In this situation $W^H = \Lie(G)^H$ is a reductive Lie algebra.
In fact, it is the Lie algebra of $Z = Z_{G}(H)$, the centralizer 
of $H$ in $G$. Note that both $N$ and $Z$ are reductive; 
cf.~\cite[Lemma 1.1]{lr}.
We claim that the natural representation 
\[ \rho \colon N \lra \GL(\Lie(Z)) \] 
is almost adjoint. Indeed, let $(H, H)$ be the commutator subgroup 
of $H$. Since $H$ acts trivially on $\Lie(Z)$,
$N^0 = Z^0 H^0$ (cf. e.g., \cite[p. 488]{lr}), 
and $(H, H) \cap Z$ is finite, we see that $\Lie(Z)$
is also the Lie algebra of $\overline{N} = N/(H,H)$ and
$\rho$ descends to the adjoint action of $\overline{N}$ 
on its Lie algebra. This proves the claim. 
Condition (i) is now established.

To check condition (ii) of Definition~\ref{def.acceptable},
we will need the following lemma.

\begin{lem} \label{lem.thm2c}
Let $G$ be a reductive (but not necessarily connected) group, 
with Lie algebra $W$ and let $G \to \GL(W)$ be the adjoint 
representation, Then 

\smallskip
{\rm(a)} the $G$-action on $W^s$ is stable for any $s \ge 1$,

\smallskip
{\rm(b)} the action of $\rho(G) = G/Z_G(G^0)$ on $W^s$ has trivial principal 
stabilizer for any $s \ge 2$.
\end{lem}

\begin{proof} (a) Since $R = \Rad(G)$ acts trivially on
$W$, the $G$-action on $W^s$ descends to an action of its
semisimple quotient $G/R$. Thus, by a theorem of
Popov~\cite{popov1} (cf. also \cite[p.~236]{pv}) it suffices
to show that $\Stab_{G/R}(w_1, \dots, w_s)$ is reductive for
$(w_1, \dots, w_s) \in W^s$ in general position.
Clearly $\Stab_{G^0/R}(w_1, \dots, w_s)$
is reductive if and only if $\Stab_{G^0}(w_1, \dots, w_s)$ is
reductive. On the other hand,
$\Stab_{G^0}(w_1, \dots, w_s)$ is contained in
$\Stab_{G^0}(w_1)$, which is a maximal torus of $G^0$,
assuming $w_1 \in W$ is in general position. Since any algebraic
subgroup of a torus is reductive, this completes the proof of part (a).

(b) By part (a) the $G$-action on $W^s$ is stable;
hence, we only need to check that the $\rho(G)$-action
on $W$ is generically free.  To do this, we may replace the adjoint
action of $G$ on its Lie algebra $W$ by the conjugation
action of $G$ on $G^0$. Indeed, applying the Luna Slice Theorem
to the 1-point orbit $\{ e \}$ in $G^0$ (with stabilizer $G$),
we see that there is a $G$-invariant open subset $U$ of $W$ containing $0$
and a $G$-equivariant morphism $U \lra G$, which is
\'etale over $e$ (and hence, dominant).  Thus if we can show
that the $G/Z_G(G^0)$-action on $(G^0)^s$ is generically free then we can
conclude that the adjoint $G$-action on $W^s$ is generically free as well.

It remains to show that a general point $(g_1, \dots, g_s) \in (G^0)^s$
has trivial stabilizer in $\rho(G) = G/Z_G(G^0)$.  In other words, we want
to show that if $g \in G$ commutes with
$g_1, \dots, g_s \in G^0$ in general position then $g$
commutes with every element of $G^0$. This
follows from~\cite[Corollary to Proposition 2]{vinberg}, which asserts that
$g_1, \dots, g_s$ generate a Zariski dense subgroup of $G^0$.
\end{proof}

We are now ready to prove that the family of representations
$\Lambda_c$ defined at the beginning of this subsection,
satisfies condition (ii) of Definition~\ref{def.acceptable}.
More precisely, we want to show that if
$\rho \colon G \to \GL(W)$ 
is an almost adjoint representat
ion and $r \ge 3$
then $W^r \catq G$ is singular away from its principal stratum.
Suppose $K \triangleleft G$ is as in Definition~\ref{def.almostadjoint}.
After replacing $G$ by $G/K$ (this doesn't change
the quotient $W^r \catq G$ or the Luna strata in it), we may
assume that $W = \Lie(G)$ and $\rho$ is the adjoint representation.

Let $\pi$ be the quotient map $W^r \to W^r \catq G$,
$x \in W^r \catq G$ be a point away from the principal stratum and
$v \in W^r$ be a point with closed $G$-orbit, such that $x = \pi(v)$.
Our goal is to show that $W^r \catq G$ is not smooth at $x$.
As in the proof of Theorem~\ref{thm2}(a) and (b) in the previous section,
we shall do this by showing that the slice representation $\Slice(v, W^r)$
is not coregular. Our strategy will be to express the $H$-representation
$\Slice(v, W^r)$
in the form \[ V_1 \oplus V_2 \oplus V_3 \, , \]
where $V_1$ has trivial principal stabilizer and $V_2$ is
not fixed pointed, then appeal to Corollary~\ref{cor1-to-prop5.1}.

Since we are assuming that $x$ does not lie in the principal stratum,
$H = \Stab_G(v) \ne \{ e \}$. As
an $H$-representation, the tangent space $T_v(G \cdot v)$ is isomorphic
to $W/\Lie(H)$ (recall that here $W = \Lie(G)$). Thus
the complement $\Slice(v, W^r)$ to $T_v(G \cdot v)$ in $W^r$
can be written as
\[ W^{r-1}  \oplus \Lie(H) \oplus S      \]
for some linear representation $S$ of $H$. By Lemma~\ref{lem.thm2c}(b),
the principal stabilizer for the $\rho(G)$-action on $W \oplus W$
is trivial. Hence, the same is true of the $H$-action on $W \oplus W$,
since $\rho(H)$ is a reductive subgroup of $\rho(G)$; see
Corollary~\ref{cor.h-m}(a).  We will now consider two cases.

\smallskip
{\bf Case 1.} $H$ acts non-trivially on $\Lie(H)$.  By
Lemma~\ref{lem.thm2c}(a), the $H$-action on $\Lie(H)$ is stable.
Since we are assuming that this action is
non-trivial, it is not fixed pointed; see
Example~\ref{ex.fixed-pointed2}.
By Corollary~\ref{cor1-to-prop5.1}, the $H$-representation
\[ \Slice(v, W^r) =
\underbrace{W^{r-1}}_{\text{\tiny trivial principal stabilizer}} \oplus
\underbrace{\Lie(H)}_{\text{\tiny not fixed pointed}} \oplus \quad S \]
is not coregular, as desired. (Recall that we are assuming throughout
that $r \ge 3$.)

\smallskip
{\bf Case 2.} $H$ acts trivially on $\Lie(H)$. Since $H$
is reductive this is only possible if $H^0$ is a central torus in $H$.
By Lemma~\ref{lem.thm2c}(b) the $G$-action (and hence, the $H$-action)
on $W^2$ has trivial principal stabilizer. Lemma~\ref{lem.central}
now tells us that the $H$-action on $W$ also has trivial principal stabilizer.
Since $H \ne \{ e \}$, no such action can be fixed pointed.
Thus by Corollary~\ref{cor1-to-prop5.1}
\[ \Slice(v, W^r) =
\underbrace{W}_{\text{\tiny trivial principal stabilizer}} \oplus
\underbrace{W}_{\text{\tiny not fixed pointed}} \oplus \quad
(W^{r-3} \oplus \Lie(H) \oplus S) \]
is not coregular.

This shows that the family $\Lambda_c$ satisfies condition (ii)
of Definition~\ref{def.acceptable}. Thus $\Lambda_c$ is an acceptable
family, and the proof of Theorem~\ref{thm2}(c) is now complete.
\qed

\section{Representation types}

Consider the action of the general linear group $\GL_n$ on the variety 
\[ V_{l, n, r} = \bbA^{lr} \times  \Mn^r \] 
% \underbrace{\Mn \times \dots \times \Mn}_{\text{\tiny $r$ times}} \, \]
given by
\[ g \cdot (a_1, \dots, a_l, A_1, \dots, A_r) \mapsto 
(a_1, \dots, a_l, gA_1g^{-1}, \dots, gA_r g^{-1}) \]
for any $a_1, \dots, a_l \in k$ and any $A_1, \dots, A_r \in \Mn$.
Let $X_{l, n, r}$ be the quotient variety $V_{l, n, r} \catq \GL_n$.
Our goal in the next two sections will be to prove Theorem~\ref{thm3} in
the following slightly more general from.

\begin{thm} \label{thm.gln}
Suppose $r \ge 3$. Then every Luna stratum 
in $X_{l, n, r} = V_{l, n, r} \catq \GL_n$ 
is intrinsic.
\end{thm}

Of course, the case where $l = 0$ is of greatest interest to us; 
in this case Theorem~\ref{thm.gln} reduces to Theorem~\ref{thm3}. 
For $l \ge 1$ the variety $X_{l, n, r}$ is only marginally 
more complicated than $X_{0, n, r}$.  Indeed, since $\GL_n$ acts 
trivially on $\bbA^{lr}$, $X_{l, n, r}$ is isomorphic to 
$\bbA^{lr} \times X_{0, n, r}$, and every Luna stratum in $X_{l, n, r}$ 
is of the form 
\begin{equation} \label{e.strata-l}
 S = \bbA^{lr} \times S_0 \, , 
\end{equation}
where $S_0$ is a Luna stratum in $X_{0, n, r}$.  We allow $l \ge 1$ 
in the statement of Theorem~\ref{thm.gln} to facilitate 
the induction step (on $n$) in the proof.

% The main result of this section, Proposition~\ref{prop.reduction3}, 
% will be a key reduction step in the proof of 
% Theorem~\ref{thm.gln} in the next section.
% In view of~\eqref{e.strata-l}, $l$ will play no role in this description.
% In particular, every assertion we make about $X_{l, n, r}$, 
% prior to Proposition~\ref{prop.reduction3} will 
% immediately reduce to the case where $l = 0$. 

The Luna stratification in $X_{l, n, r}$, has a natural 
combinatorial interpretation, which we will now recall, in preparation 
for the proof of Theorem~\ref{thm.gln} in the next section. 
To $x \in X_{l, n, r}$, one can uniquely associate
a point $v = (a, \dots, a_l, A_1, \dots, A_r) \in V_{l, n, r}$,
with a closed $\GL_n$-orbit.  Here each $a_i \in k$ and each 
$A_j \in \Mat_r$.  We will view an $r$-tuple $(A_1, \dots, A_r) \in \Mn^r$ 
of $n \times n$-matrices as an $n$-dimensional representation 
\[ \rho \colon k \{ x_1, \dots, x_r \} \to \Mn \]
of the free associative algebra $k \{ x_1, \dots, x_r \}$ 
on $r$ generators, sending $x_i$ to $A_i$. By a theorem 
of Artin~\cite[(12.6)]{artin}, the orbit of $v$
is closed in $\Mn^r$ if and only if $\phi$ is semisimple. 
(Strictly speaking, Artin's theorem only covers the case where $l = 0$;
but since $V_{l, n, r} = \bbA^{lr} \times V_{0, n, r}$ and $\GL_n$ acts
trivially on $\bbA^{lr}$, the general case is an immediate consequence.)
If $\rho$ can be written as $\rho_1^{e_1} \oplus \dots \oplus \rho_s^{e_s}$,
where 
\[ \rho_i \colon k \{ x_1, \dots, x_r \} \to \M_{d_i} \]
is an irreducible $d_i$-dimensional representation, we will say
that the {\em representation type} of $x$ is 
\begin{equation} \label{e.rep-type}
\tau = [(d_1, e_1), \dots, (d_r, e_r)] \, .
\end{equation}
The square brackets $[ \quad ]$ are meant to indicate that
that $\tau$ is an {\em unordered} collection of pairs $(d_i, e_i)$;
permuting these pairs does not change the representation type.
Note also 
that $d_i, e_i \ge 1$ for every $i = 1, \dots, s$.  Following le Bruyn 
and Procesi~\cite{lbp}, we shall denote the set of representation 
types~\eqref{e.rep-type} with $d_1e_1 + \dots + d_s e_s = n$ by $\RT_n$. 
If $\tau = [(d_1, e_1), \dots, (d_s, e_s)] \in RT_n$ and 
$\mu =  [(d'_1, e'_1), \dots, (d'_{s'}, e'_{s'})] \in \RT_{n'}$ then
we will sometimes denote the representation type
\[ [(d_1, e_1), \dots, (d_s, e_s),  
(d'_1, e'_1), \dots, (d'_{s'}, e'_{s'})] \in \RT_{n + n'} \]
by $[\tau, \mu]$.

The Luna strata in $X_{l, n, r} = (\bbA^{lr} \times \Mn^r) \catq \GL_n$ 
are in a 1-1 correspondence 
with $\RT_n$; cf.i, e.g.,~\cite[Section 2]{lbp}.
If $x \in X_{l, n, r}$ 
has representation type $\tau$, as in~\eqref{e.rep-type}, 
then the associated stabilizer subgroup 
\begin{equation} \label{e.stab_tau}
H_{\tau} = \Stab_G(v) \simeq \GL_{e_1} \times \dots \times \GL_{e_s} \, , 
\end{equation}
embedded into $\GL_n$ as follows.  Write
\[ k^n = V_1 \otimes W_1\oplus \dots \oplus V_s \otimes W_s \, , \]
where $\dim(V_i) = d_i$ and $\dim(W_i) = e_i$ and let $\GL_{e_i}$ 
act on $W_i$. Then 
$H_{\tau} = \GL_{e_1} \times \dots \times \GL_{e_s}$ 
is embedded in $\GL_n$ via 
\[ (g_1, \dots, g_s) \mapsto 
(I_{d_1} \otimes g_1 \oplus \dots \oplus I_{d_s} \otimes g_s) \, . \] 
where $I_d$ denotes the $d \times d$ identity matrix; 
cf.~\cite[Section 2]{lbp}.  For notational convenience 
we shall denote the Luna strata in $X_{l, n, r} = \Mn^r \catq \GL_n$ by 
$X_{n, r}^{\tau}$, rather than $X_{n, r}^{(H_{\tau})}$. Note that if 
$\tau = [(d_1, e_1), \dots, (d_s, e_s)] \in \RT_n$ then
\begin{equation} \label{e.dimension}
\dim \, X_{l, n, r}^{\tau} = (r-1)(d_1^2 + \dots + d_s^2) + s + lr 
\end{equation}
for any $r \ge 2$; cf.~\cite[p. 158]{lbp}.

\begin{defn} \label{def.partial-order}
An {\em elementary refinement} of 
$\tau = [(d_1, e_1), \dots, (d_s, e_s)] \in \RT_n$ consists in either

\smallskip
(1) replacing the pair $(d_i, e_i)$ by two pairs $(a_i, e_i)$ 
and $(b_i, e_i)$, where $a_i, b_i \ge 1$ and $a_i + b_i = d_i$ or

\smallskip
(2) replacing two pairs $(d_i, e_i)$ and $(d_j, e_j)$, 
with $d_i = d_j$, by the single pair $(d_i, e_i + e_j)$.

\smallskip
\noindent
Given two representation types $\tau$ and $\tau'$, we will say that
$\tau' \prec \tau$ if $\tau'$ can be obtained from $\tau$ by 
a sequence of elementary refinements. This defines a partial order
$\preceq$ on $\RT_n$.
\end{defn}

Note that while operations (1) and (2) are defined in purely 
combinatorial terms, they are, informally speaking, 
designed to reflect the two ways a representation
\[ \rho = \rho_1^{e_1} \dots \rho_s^{e_s} \colon
k \{ x_1, \dots, x_r \} \to \Mn \]
can ``degenerate". Here $\rho_1, \dots, \rho_s$ are distinct
irreducible representations of dimensions $d_1, \dots, d_s$ respectively.
In case (1), one of the representations $\rho_i$ ``degenerates" 
into a direct sum of irreducible subrepresentations of degree
$a_i$ and $b_i$ (each with multiplicity $e_i$). In case (2) 
we ``degenerate" $\rho$ by making $\rho_i$ and $\rho_j$ isomorphic
(of course, this is only possible if their dimensions $d_i$ and $d_j$ are
the same). The following lemma gives this a precise meaning.

\begin{lem} \label{lem.order} $X_{n, r}^{\mu}$ lies 
in the closure of $X_n^{\tau}$ if and only if  $\mu \preceq \tau$. 
\end{lem}

\begin{proof} In view of \eqref{e.strata-l}, we may assume $l = 0$. 
In this case Lemma~\ref{lem.order} is proved 
in~\cite[Theorem II.1.1]{lbp}.
\end{proof}
%%%%%%%%%%%%%%%%%%%%%%%%%%%%%%%%%% 
% In particular, 
% \begin{equation} \label{e.closure}
% \overline{X_{l, n, r}^{\tau}} = 
% \bigcup_{\tau' \preceq \tau} X_{l, n, r}^{\tau'} \, . 
% \end{equation}
%%%%%%%%%%%%%%%%%%%%%%%%%%%%%%%%%% 

\begin{remark} \label{rem.graph}
Note that $V_{l, n, r} = W^r$, where
$W = \bbA^l \times \Mn$ is the Lie algebra 
of $G = (\GL_1)^l \times \GL_n$.  The $(\GL_1)^l$ factor acts trivially 
on $W$ (via the adjoint action), so we may drop 
it without changing the quotient $W^r \catq G$ or the 
Luna strata in it. Theorem~\ref{thm2}(c) now tells us that
the Luna stratification in
$X_{l, n, r} = V_{l, n, r} \catq \GL_n = W^r \catq G$ is intrinsic
(provided
$r \ge 3$). Thus an automorphism $\sigma$ of $X_{l, n, r}$
induces an automorphism $\sigma^*$ of the set $\RT_n$ 
of Luna strata in $X_{l, n, r}$ given 
by $\sigma^*(\tau) = \nu$ if $\sigma(X^{\tau}) = X^{\nu}$.
By Lemma~\ref{lem.order} $\sigma^*$ respects the partial 
order $\preceq$ on $\RT_n$. Theorem~\ref{thm.gln} asserts 
that $\sigma^*$ is always trivial. In some cases this can be 
deduced from the fact that the partially ordered set 
$(\RT_n, \preceq)$ has no non-trivial automorphisms.
For example, the partially ordered sets $\RT_n$, for $n = 1$, $2$ 
and $3$, pictured below have no non-trivial automorphisms. 
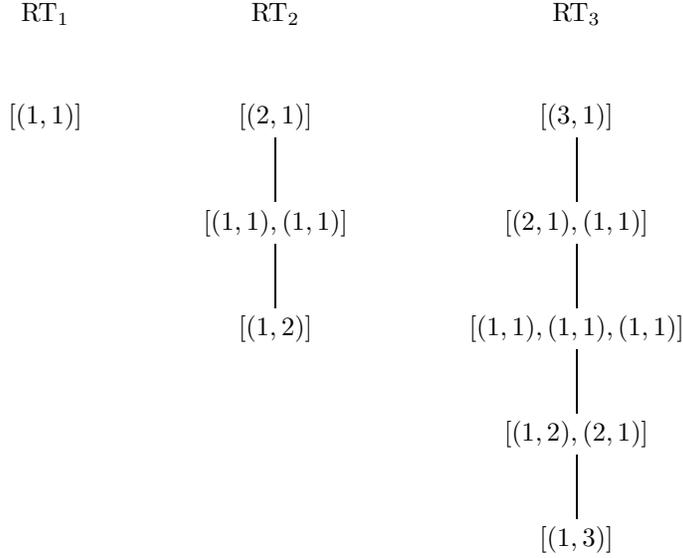
\[ \xymatrix{ \RT_1 \cr
            [(1, 1)]} 
\quad \quad
\quad \quad
 \xymatrix{ \RT_2 \cr
            [(2, 1)] \ar@{-}[d] \cr 
            [(1, 1), (1, 1)] \ar@{-}[d] \cr 
            [(1, 2)]}
\quad \quad
\quad \quad
 \xymatrix{ \RT_3 \cr
            [(3, 1)] \ar@{-}[d] \cr 
            [(2, 1), (1, 1)] \ar@{-}[d] \cr 
            [(1, 1), (1, 1), (1, 1)] \ar@{-}[d] \cr 
            [(1, 2), (2, 1)] \ar@{-}[d] \cr 
            [(1, 3)]} \]
This proves Theorem~\ref{thm.gln} for $n \le 3$.
\end{remark}

Note however, that for larger $n$ the partially ordered set 
$(\RT_n, \preceq)$ does have non-trivial automorphisms. For example,
a quick look at $\RT_4$ (pictured on p. 156 in~\cite{lbp}), shows that
the permutation $\alpha$ of $\RT_4$ interchanging 
\[ \text{$\tau = [(1, 1), (1, 1), (1, 1), (1, 1)]$
and $\nu = [(2, 1), (1, 2)]$} \]
and fixing every other element, does, indeed, respect the partial order 
on $\RT_4$. For this reason we cannot hope to prove Theorem~\ref{thm.gln}
by purely combinatorial arguments, without taking into account 
the geometry of the strata $X_{l, n, r}^{\tau}$. 
%%%%%%%%%%%%%%%%%%%%%%%%%%%%%%%%%%%%%%%%
% (In the above-mentioned example, we note
% that the strata have different dimensions, $\dim(X^{\tau}) = 4(r-1) + 4 + lr$
% but $\dim(X^{\tau}) = 5(r-1) + 2 + lr$
%%%%%%%%%%%%%%%%%%%%%%%%%%%%%%%%%%%%%%%
Nevertheless, the following combinatorial proposition will play a key
role in the proof of Theorem~\ref{thm.gln} in the next section.

\begin{prop} \label{prop.reduction3} 
{\rm(a)} Let $\alpha$ be an automorphism of $\RT_n$ {\rm(}as a partially ordered set{\rm)}.
If $\alpha([(1, 1), \mu]) = [(1, 1), \mu]$
for every $\mu \in \RT_{n-1}$ then $\alpha = \id$.

\smallskip
{\rm(b)} Let $l \ge 0$, $n \ge 1$ and $r \ge 3$. Suppose we know that Luna strata
of the form $X_{l, n, r}^{[(1, 1), \mu]}$ are intrinsic in $X_{l, n, r}$
for any $\mu \in RT_{n-1}$.  Then every Luna stratum in $X_{l, n, r}$ 
is intrinsic.
\end{prop}

\begin{proof} Part (b) is an immediate consequence of part (a) and
Remark~\ref{rem.graph}; we shall thus concentrate on proving part (a).
Given a representation type
$\tau = [(d_1, e_1), \dots, (d_s, e_s)]$, let $m(\tau)$
denote the minimal value of $d_i + e_i$, as $i$ ranges from $1$ to $s$.
Note that since $d_i, e_i \ge 1$ for each $i$, we have $m(\tau) \ge 2$.
We will now show that $\alpha(\tau) = \tau$ by induction on $m(\tau)$. 
By our assumption this is the case if $m(\tau) = 2$, since in this 
case $d_i = e_i = 1$ for some $i$.

For the induction step, assume that $m(\tau) = m \ge 3$ and
$\alpha(\nu) = \nu$ for every $\nu \in \RT_n$ with $m(\nu) < m$.
After renumbering
the pairs $(d_i, e_i)$, we may assume that $d_1 + e_1 = m$.
Suppose 
\begin{equation} \label{e.tau'} 
\alpha(\tau) = \tau' = [(d_1', e_1'), \dots, (d_t', e_t')] \, . 
\end{equation}
Our goal is to show that $\tau = \tau'$. We will consider 
two cases, where $d_1 \ge 2$ and $e_1 \ge 2$ separately.
Since we are assuming that $d_1 + e_1 = m \ge 3$, 
these two cases cover every possibility.

\smallskip
{\bf Case 1:} $d_1 \ge 2$.  Let 
\[ \tau_0 = [(1, e_1), (d_1-1, e_1), (d_2, e_2), \dots, (d_s, e_s)] \, . \]
Since $m(\tau_0) = e_1 + 1 < d_1 + e_1 = m$, the induction assumption 
tells us that $\alpha(\tau_0) = \tau_0$. 
Now observe that $\tau_0$ immediately precedes $\tau$ in the partial 
order on $\RT_n$, i.e., $\tau_0$ is obtained from 
$\tau$ by a single elementary refinement;  
see Definition~\ref{def.partial-order}. 
Consequently, $\tau_0$ can also be obtained from $\tau'$ by a single 
elementary refinement. Schematically,
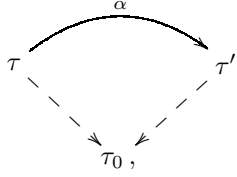
\[ \xymatrix{ \tau \ar@{-->}[dr]  \ar@/^1.5pc/[rr]^{\alpha} &  
& \tau' \ar@{-->}[dl] \cr 
 & \tau_0 \, ,&  } \]
where the broken arrows denote elementary refinements.
Now observe that 
\[ (d_i', e_i') \neq (1, e_1) \; \; \; \text{or} \; \; \; (d_1 -1, e_1) \quad 
\forall \quad i = 1, \dots, s \, . \]
Indeed, otherwise we would have $m(\tau') < m$ and
thus $\alpha(\tau') = \tau'$ by the induction assumption.
Combining this with~\eqref{e.tau'}, we obtain
$\alpha(\tau) = \alpha(\tau') = \tau'$. Since $\alpha$ 
is a permutation of $\RT_n$, we conclude that $\tau = \tau'$, contradicting
$m(\tau') < m$.
%  (Recall that $\alpha$
% is assumed to be a permutation of $\RT_n$.) 

To sum up, $\tau_0$ ``contains" two pairs, 
$(1, e_1)$ and $(d_1-1, e_1)$, that are not ``present" in $\tau'$.
It now follows from Definition~\ref{def.partial-order} that the only possible
elementary refinement taking $\tau'$ to $\tau_0$ is of type (1),
consisting of "splitting up" $(d_1, 1)$ into $(d_1 -1, 1)$ and $(1, 1)$. 
That is, $\tau' = [(d_1, 1), \mu] = \tau$, as claimed.

\smallskip
{\bf Case 2:} $e_1 \ge 2$. The argument here is very similar (or more 
precisely, ``dual", in the sense explained in Remark~\ref{rem.duality}) 
to the one in Case 1. Let
\[ \tau_1 = [(d_1, 1), (d_1, e_1 -1), (d_2, e_2), \dots, (d_s, e_s)] \, . \]
Since $m(\tau_1) = d_1 + 1 < d_1 + e_1 = m$, the induction assumption 
tells us that $\alpha(\tau_1) = \tau_1$. 
The relationship between $\tau$, $\tau'$ 
and $\tau_1$ is shown in the following diagram
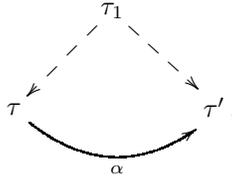
\[ \xymatrix{ & \tau_1 \ar@{-->}[dr] \ar@{-->}[dl] &  \cr
\tau  \ar@/_1.5pc/[rr]_{\alpha} &  & \tau' \, ,} \]
where the broken arrows denote elementary refinements.  Once again, 
we see that $\tau_1$ ``contains" two pairs, 
$(d_1, 1)$ and $(d_1, e_1 - 1)$ both of which "disappear" after we
perform an elementary refinement (and obtain $\tau'$).
% \[ (d_i', e_i') \neq (d_1, 1) \; \; \; \text{or} \; \; \;  
% (d_1, e_1 - 1) \quad \forall \quad i = 1, \dots, s \, \]
% since otherwise we would have $m(\tau') < m$ and
% thus, by the induction assumption, $\alpha(\tau') = \tau'$. 
% Combining this with~\eqref{e.tau'}, we 
% obtain $\tau = \tau'$, a contradiction (since we are assuming 
% that $m(\tau') < m$).
% Thus $\tau_1$ ``contains" two pairs, 
% $(d_1, 1)$ and $(d_1, e_1 - 1)$ both of which "disappear" after we
% perform an elementary refinement (and obtain $\tau'$).
This is only possible if the elementary refinement 
taking $\tau_1$ to $\tau'$ is of type (2) 
(cf. Definition~\ref{def.partial-order}) 
and consists of replacing 
$(d_1, 1)$ and $(d_1, e_1 - 1)$ by $(d_1, e_1)$.
This shows that $\tau = \tau'$, thus
completing the proof of Proposition~\ref{prop.reduction3}.
\end{proof}

\begin{remark} \label{rem.duality}
The two elementary refinement operations 
of Definition~\ref{def.partial-order} are dual
to each other in the following sense. Given a representation type
\[ \tau = [(d_1, e_1), \dots, (d_s, e_s)] \]
% we will use the 
% symbol $\overline{\tau}$ for the representation type
let $\overline{\tau} = [(e_1, d_1), \dots, (e_s, d_s)]$. 
% and
% the symbol $\tau \xrightarrow{(i)} \nu$ to indicate that $\nu$ is 
% obtained from $\tau$ by a single elementary refinement of type ($i$),
% where $i = 1$ or $2$. Using this terminology, duality between 
% the elementary refinement 
% operations (1) and (2) can be expressed as follows:
% \[ \text{$\tau \xrightarrow{(1)} \nu$ if and only if $\overline{\nu} 
% \xrightarrow{(2)} \overline{\tau}$ 
% and $\tau \xrightarrow{(2)} \nu$ if and only if 
% $\overline{\nu} \xrightarrow{(1)} \overline{\tau}$.} \]
% \[ \xymatrix{ \tau \ar@{->}[d]^{(1)}  \cr \nu \, ,}$ 
% \quad \text{if and only if} \quad
% if and only if
% $\xymatrix{ \overline{\nu} \ar@{->}[d]^{(2)} \cr \overline{\tau} \, ,}$ 
% % \quad \text{and} \quad \text{if} 
% and 
% $\xymatrix{ \tau \ar@{->}[d]^{(2)}  \cr \nu \, ,}$ 
% % \quad \text{then} \qThe geometric meaning of this duality, i.e., the
relationship between the strata
$X_{l, n, r}^{\tau}$ and $X_{l, n, r}^{\overline{\tau}}$, is not
entirely clear to us. 
% Note, in particular, that these strata 
% are usually of different dimension; cf.~\eqref{e.dimension}. uad
% if and only if
% $\xymatrix{ \overline{\nu} \ar@{->}[d]^{(1)} \cr \overline{\tau} \, ,}$ 
Then $\alpha$ is 
 obtained from $\beta$ by an elementary refinement of type (1) 
 (respectively, of type (2)) if and only if $\overline{\beta}$ is obtained 
 from $\overline{\alpha}$ by an elementary refinement of type (2) 
 (respectively, of type (1)). 
Consequently, the map $\tau \mapsto \overline{\tau}$
is an isomorphism between the partially ordered sets $(\RT_n, \succeq)$
and $(\RT_n, \preceq)$. The statement and the proof of  
Proposition~\ref{prop.reduction3} are invariant with respect 
to this map (in particular, Case 2 is dual to Case 1).
\end{remark}

\section{Proof of Theorem~\ref{thm3}}

In this section we will prove Theorem~\ref{thm.gln}. This will immediately
yield Theorem~\ref{thm3} (for $l = 0$). We will continue to use
the notations introduced in the previous section; in particular,
$V_{l, n, r}$ stands for $\bbA^{lr} \times \Mn$ and 
$X_{l, n, r}$ denotes the categorical quotient $V_{l, n, r} \catq \GL_n$.

We will argue by induction on $n$. The base cases, $n = 1$ and $2$, 
are proved in Remark~\ref{rem.graph}. (Theorem~\ref{thm.gln} 
is also proved there for $n = 3$ but we shall not need that here.)
For the induction step assume $n \ge 3$ and
$\sigma$ is an automorphism of $X_{l, n, r}$.
Recall that $\sigma$ maps each stratum in $X_{l, n, r}$ 
to another stratum; cf. Remark~\ref{rem.graph}. In particular, $\sigma$ 
preserves the maximal (principal) stratum 
$X_{l, n, r}^{[(n, 1)]}$ (which is the unique stratum 
of maximal dimension) and permutes the "submaximal" strata  
$X_{l, n, r}^{[(d, 1), (n-d, 1)]}$, $1 \le d \le \frac{n}{2}$, 
among themselves.  By the dimension formula~\eqref{e.dimension},
\[ \dim \, X^{[(d, 1), (n-d, 1)]} = rl + 2 + (r-1)(d^2 + (n-d)^2) 
= rl + 2 + 2(r-1)((d - \frac{n}{2})^2 + \frac{n^2}{4}) \, . \]
Thus the submaximal strata 
$X_{l, n, r}^{[(d, 1), (n-d, 1)]}$ have different dimensions 
for different values of $d$ between $1$ and $\frac{n}{2}$. 
Hence, $\sigma$ preserves each one of them.

Of particular interest to us is the submaximal stratum 
$X_{l, n, r}^{[(1, 1), (n-1, 1)]}$.  Since $\sigma$ 
preserves this stratum, it preserves its closure
$\overline{X_{l, n, r}^{[(1, 1), (n-1, 1)]}}$ 
and thus lifts to an automorphism $\tilde{\sigma}$ of
the normalization of $\overline{X_{l, n, r}^{[(1, 1), (n-1, 1)]}}$. 
The rest of the argument will proceed as follows.  We will identify 
the normalization of $\overline{X_{l, n, r}^{[(1, 1), (n-1, 1)]}}$ 
with $X_{l+1, n-1, r}$ and 
relate Luna strata in $X_{l+1, n-1, r}$ and $X_{l, n, r}$ via
the normalization map.
By the induction assumption $\tilde{\sigma}$ preserves 
every Luna stratum in $X_{l+1, n-1, r}$; using the normalization map
we will be able to conclude that $\sigma$ preserves certain Luna strata 
in $X_{l. n, r}$. Proposition~\ref{prop.reduction3} will then
tell us that, in fact, $\sigma$ preserves 
every Luna stratum in $X_{l. n, r}$, thus completing the proof.

We now proceed to fill in the details of this outline. 
First will explicitly describe the normalization map
for $\overline{X_{l, n, r}^{[(1, 1), (n-1, 1)]}}$.  The stabilizer
$H_{\tau}$, corresponding to $\tau = [(1, 1), (n-1, 1)]$ 
consists of matrices of the form 
\[ \diag(a, \underbrace{b, \dots,b}_{\text{\tiny $n-1$ times}}) \, , \]
cf.~\eqref{e.stab_tau}. General theory now tells us that 
the natural projection 
\[ V_{l, n, r}^{H_{\tau}} \catq N_G(H_{\tau}) 
\to \overline{X_{l, n, r}^{\tau}} \] 
is the normalization map for $\overline{X_{l, n, r}^{\tau}}$. Here  
\[ V_{l, n, r}^{H_{\tau}} = \bbA^{lr} \times 
  \Mat_1^r \times \Mat_{n - 1}^r = \bbA^{l(r+1)} \times \Mat_{n-1}^r = 
V_{l + r, n-1, r} \] 
and 
\begin{equation} \label{e.N_tau}
 N_G(H_{\tau}) = \GL_1 \times \GL_{n-1} \, , 
\end{equation}
where $\GL_1$ acts trivially on $\bbA^{l(r + 1)} \times \Mat_{n-1}^r$ 
and $\GL_{n-1}$ acts on the second factor by simultaneous conjugation.
Since $\GL_1$ acts trivially, we may replace $\GL_1 \times \GL_{n-1}$
by $\GL_{n-1}$ without changing the categorical quotient or the Luna 
strata in it. 
That is, the normalization $V_{l, n, r}^{H_{\tau}} \catq N_G(H_{\tau})$
of $\overline{X_{l, n, r}^{[(1, 1), (n-1, 1)]}}$ is (canonically) 
isomorphic to $V_{l+1, n-1, r} \catq \GL_{n-1} = X_{l+1, n-1, r}$. 

The following lemma gives a summary of this construction. Here,
as before, we identify an $r$-tuple $A = (A_1, \dots, A_r)$
of $n \times n$-matrices with the $n$-dimensional representation
$\rho_A \colon k \{ x_1, \dots, x_r \} \to \Mn$ of the free
associative $k$-algebra $k \{ x_1, \dots, x_n \}$, taking $x_i$ to $A_i$. 

\begin{lem} \label{lem.norm}
 Suppose $n \ge 3$ and $r \ge 2$ and let
 $f \colon V_{l+1, n-1, r} \to V_{l, n, r}$
 be the morphism given by
 \[ f \colon (t_1, \dots, t_{lr}, \rho) \mapsto
 (t_1, \dots, t_{l(r-1)}, \rho_t \oplus \rho)) \, , \]
 where $t = (t_{l(r-1)+1}, t_{l(r-1) + 2}, \dots, t_{lr}) \in \Mat_1^r$.
 Then

 \smallskip
 {\rm(a)} $f$ descends to the normalization map
 \[ \overline{f} \colon X_{l+1, n-1, r} \to
 \overline{X_{l, n, r}^{[(1, 1), (n-1, 1)]}} \, , \]
where $X_{l, n, r} = V_{l, n, r} \catq \GL_n$.

 {\rm(b)} $\overline{f}$ maps $\overline{X_{l+1, n-1, r}^{\mu}}$
 onto $\overline{X_{l, n, r}^{[(1, 1), \mu]}}$
 for every $\mu \in \RT_{n-1}$.
\end{lem}

\begin{proof} Part (a) follows from the discussion before 
the statement of the lemma. To prove part (b), observe 
that every semisimple representation 
$k \{ x_1, \dots, x_r \} \to \Mn$ of 
type $[(1, 1), \mu]$ can be written in the form
$\rho_0 \oplus \rho$, where
$\rho$ is an $n-1$-dimensional representation of type $\mu$ and $\rho_0$ 
is a 1-dimensional representation (of type $(1, 1)$).
This shows that the image of $X_{l+1, n-1, r}^{\mu}$ contains 
$X_{l, n, r}^{[(1, 1), \mu]}$. Since
\[ \dim \, X_{l+1, n-1, r}^{\mu} = \dim
X_{l, n, r}^{[(1, 1), \mu]} \, , \]
see \eqref{e.dimension}, and $\overline{f}$ is a finite map (in particular, 
$\overline{f}$ takes closed sets to closed sets), part (b) follows.
\end{proof}

\begin{remark} Lemma~\ref{lem.norm} uses, in a crucial way, 
the assumption that
$n \ge 3$.  If $n = 2$ then~\eqref{e.N_tau} fails; instead we have 
$N_G(H_{\tau}) = (\GL_1 \times \GL_1) \sdp \Sym_2$, and the 
entire argument falls apart.
\end{remark}

We are now ready to finish the proof of Theorem~\ref{thm.gln}. Restricting
$\sigma$ to the closure of the stratum $X_{l, n, r}^{[(1, 1), (n-1, 1)]}$
and lifting it an automorphism $\tilde{\sigma}$
to the normalization, we obtain the following 
commutative diagram:
\[ \xymatrix{ X_{l+1, n-1, r} \ar@{->}[r]^{\tilde{\sigma}} 
\ar@{->}[d]_(.4){\overline{f}} & 
X_{l+1, n-1, r} \ar@{->}[d]^(.4){\overline{f}} \cr 
 \overline{X_{l, n, r}^{[(1, 1), (n-1, 1)]}} 
 \ar@{->}[r]^{\sigma} &
 \overline{X_{l, n, r}^{[(1, 1), (n-1, 1)]}} \, .} \]
By our induction assumption $\tilde{\sigma}$ preserves every
Luna stratum $X_{l+1, n-1, r}^{\mu}$ in $X_{l+1, n-1, r}$.
Hence, by Lemma~\ref{lem.norm}(b), $\sigma$ preserves the closure
of every Luna stratum in $X_{l, n, r}$ of the form
$X_{l, n, r}^{[(1, 1), \mu]}$. Since 
$X_{l, n, r}^{[(1, 1), \mu]}$ is the unique stratum of maximal dimension in
its closure, we conclude that $\sigma$ preserves
$X_{l, n, r}^{[(1, 1), \mu]}$ for every $\mu \in \RT_{n-1}$.
By Propostion~\ref{prop.reduction3}(b), we conclude that 
$\sigma$ preserves every Luna stratum in $X_{l, n, r}$. 
This concludes the proof of Theorem~\ref{thm.gln} and thus of
Theorem~\ref{thm3}.
\qed

\begin{remark} \label{rem.(2,2)} 
Theorem~\ref{thm3} fails if (a) $r = 1$ or (b) $(n, r) = (2, 2)$,
because in this case the ring of invariants $R = k[\Mn^r]^{\GL_n}$
is a polynomial ring or equivalently, $X = \Mn^r \catq \GL_n$ is
an affine space. In case (a) $R$ is freely generated by 
the coefficients of the characteristic polynomial of $A \in \M_n$,
viewed as $\GL_n$-invariant polynomials $\Mn^1 \to k$ and
in case (b) by the five $\GL_2$-invariants $\Mn^2 \to k$ given by
$(A_1, A_2) \mapsto \tr(A_1)$, $\tr(A_2)$, $\det(A_1)$, $\det(A_2)$, 
and $\det(A_1 + A_2)$, respectively; cf. \cite[VIII, Section 136]{gy}.  
\end{remark}

\section*{Acknowledgements} We would like to thank J.-L. Colliot-Th\'el\`ene,
Ph. Gille and M. Thaddeus for helpful discussions. We are also grateful 
to the Pacific Institute for the Mathematical Sciences (PIMS) for 
making this collaboration possible.

\end{document}